\def\IS{{\Bbb S}} 
 
\def\IK{{\Bbb K}}

\def\IC{\Bbb C} 
\def\ID{{\Bbb D}}

\def\zbar{{\overline{z}}} 
 
\def\wbar{{\overline{w}}} 
\def\mubar{{\overline{\mu}}}

\documentclass[10pt]{article} 

\usepackage[psamsfonts]{amssymb} 
\usepackage{amsfonts,amsmath}  

\usepackage{epsfig,multicol}

\newtheorem{theorem}{Theorem}[section] 
 
\newtheorem{lemma}{Lemma}[section] 
\newtheorem{corollary}{Corollary}[section]

\newtheorem{conjecture}{Conjecture}[section]

\numberwithin{equation}{section}

\title{{\bf The $L^p$ Teichm\"uller theory:}\\ { Existence and regularity of critical points}}
\author{Gaven Martin \& Cong Yao \thanks{ This research of both authors is supported in part by a grant from the NZ Marsden Fund.\newline
Part of this work appears in the second author's PhD Thesis.
\newline
{\bf  Mathematics subject classification 2010,}{Primary: 30C62 31A05 49J10;  } \newline 
{\bf  Keywords,} {Calculus of variations,  quasiconformal, distributional equations, mean distortion}}}
\date{}
\begin{document}
\maketitle

\begin{abstract}  
We study minimisers of the $p$-conformal energy functionals,
\[
\mathsf{E}_p(f):=\int_\ID \IK^p(z,f)\,dz,\quad f|_\IS=f_0|_\IS,
\]
defined for self mappings $f:\ID\to\ID$ with finite distortion and prescribed boundary values $f_0$.  
Here 
\[  \IK(z,f) = \frac{\|Df(z)\|^2}{J(z,f)} = \frac{1+|\mu_f(z)|^2}{1-|\mu_f(z)|^2}\]  
is the pointwise distortion functional and $\mu_f(z)$ is the Beltrami coefficient of $f$. We show that for quasisymmetric boundary data the limiting regimes $p\to\infty$ recover the classical Teichm\"uller theory of extremal quasiconformal mappings (in part a result of Ahlfors), and for $p\to1$ recovers the harmonic mapping  theory.

Critical points of $\mathsf{E}_p$ always satisfy the inner-variational distributional equation 
\[
2p\int_\ID \IK^p\;\frac{\overline{\mu_f}}{1+|\mu_f|^2}\varphi_\zbar \; dz=\int_\ID \IK^p \; \varphi_z\; dz,\quad\forall\varphi\in C_0^\infty(\ID ).
\]
We establish the existence of minimisers in the {\em a priori} regularity class $W^{1,\frac{2p}{p+1}}(\ID)$ and show these minimisers have a pseudo-inverse - a continuous $W^{1,2}(\ID)$ surjection of $\ID$ with $(h\circ f)(z)=z$ almost everywhere.  We then give a sufficient condition to ensure $C^{\infty}(\ID)$ smoothness of solutions to the distributional equation.  For instance $\IK(z,f)\in L^r_{loc}(\ID)$ for any $r>p+1$ is enough to imply the solutions to the distributional equation are local diffeomorphisms. Further $\IK(w,h)\in L^1(\ID)$ will imply $h$ is a homeomorphism,  and together these results yield a diffeomorphic minimiser.  We show such higher regularity assumptions to be necessary for critical points of the inner variational equation.  \end{abstract}

\section{Introduction}

A mapping $f:\ID\to \ID$ has finite distortion if 
\begin{enumerate}
\item $f\in W^{1,1}_{loc}(\ID)$,  the Sobolev space of functions with locally integrable first derivatives,
\item the Jacobian determinant $J(z,f)\in L^{1}_{loc}(\ID)$, and 
\item there is a measurable function ${\bf K}(z)\geq 1$, finite almost everywhere, such that 
 \begin{equation}\label{1.1}
 |Df(z)|^2 \leq {\bf K}(z) \, J(z,f), \hskip10pt \mbox{ almost everywhere in $\ID$}.
 \end{equation}
\end{enumerate}
See \cite[Chapter 20]{AIM} for the basic theory of mappings of finite distortion and the associated governing equations; degenerate elliptic Beltrami systems.
In (\ref{1.1}) the operator norm is used.  However this norm loses smoothness at crossings of the eigenvalues and for this reason when considering minimisers of distortion functionals one considers the distortion functional
\begin{equation}\label{1.2}
\IK(z,f) = \left\{\begin{array}{cc} \frac{\|Df(z)\|^2}{J(z,f)}, & \mbox{if $J(z,f)\neq 0$} \\
1, & \mbox{if $J(z,f)= 0$.} \end{array}\right.
\end{equation}
This was already realised by Ahlfors in his seminal work proving Teichm\"uller's theorem and establishing the basics of the theory of quasiconformal mappings, \cite[\S 3, pg 44]{Ahlfors}. We reconcile (\ref{1.1}) and (\ref{1.2}) by noting $\IK(z,f)=\frac{1}{2}\big({\bf K}(z)+1/{\bf K}(z) \big)$ almost everywhere,  where ${\bf K}(z)$ is chosen to be the smallest functions such that (\ref{1.1}) holds.

\subsection{Minimising mean distortion.}

Let $p\geq1$.  The $L^p$ mean distortion of a self-homeomorphism of $\overline{\ID}$ and of finite distortion is defined as
\begin{equation}\label{1.3}
\mathsf{E}_p(f):=\int_\ID \IK^p(z,f)\; dz,
\end{equation}
Note
\begin{equation}\label{1.4}
\IK(z,f)=\frac{1+|\mu_f|^2}{1-|\mu_f|^2}=\frac{|f_z|^2+|f_\zbar|^2}{|f_z|^2-|f_\zbar|^2}
\end{equation}
where $\mu_f=f_\zbar/f_z$ is the Beltrami coefficient of $f$.  This follows from (\ref{1.1}).  For brevity we write $\IK_f=\IK(z,f)$.

\medskip

 Let $f_0:\overline{\ID }\to\overline{\ID }$ be a finite distortion homeomorphism such that $\mathsf{E}_p(f_0)<\infty$. We regard $f_0$ as the boundary data and define the space of functions
\[ 
\mathsf{F}_p:=\Big\{f\in W_{loc}^{1,1}(\ID ):\mathsf{E}_p(f)<\infty,\; f|_\IS=f_0|_\IS, \; \mbox{$f:\overline{\ID }\to \overline{\ID}$  is a homeomorphism} \Big\}
\]
Every mapping in $\mathsf{F}_p$ has finite distortion.  We recall the following conjecture announced in \cite{IMO}.
\begin{conjecture}\label{conjecture}
In the space $\mathsf{F}_p$, there is a minimiser $f$ such that
\[
\mathsf{E}_p(f)=\min_{g\in\mathsf{F}_p}\mathsf{E}_p(g).
\]
Furthermore, this map is a $C^\infty$-smooth diffeomorphism from $\ID$ to $\ID$.
\end{conjecture}
This conjecture contains two parts. First, in the space $\mathsf{F}_p$ there is a minimiser; and second, if there is a homeomorphic minimiser, it must be a diffeomorphism.  There is some evidence to support this conjecture in \cite{MJ} where it is shown that the extremals for the similar $L^p$-Gr\"otzsch problem are unique and are diffeomorphisms. However our examples below point in the other direction.

\medskip

It is proved in \cite{IMO} that a minimiser in $\mathsf{F}_p$ must satisfy the following inner-variational equation:
\begin{equation}\label{1.5}
2p\int_\ID \IK_f^p\frac{\overline{\mu_f}}{1+|\mu_f|^2}\varphi_\zbar dz=\int_\ID \IK_f^p\varphi_zdz,\quad\forall\varphi\in C_0^\infty(\ID ).
\end{equation}
This equation arises as follows.  Let $\varphi\in C^{\infty}_{0}(\ID)$ with $\|\nabla\varphi\|_{L^{\infty}(\ID)}<1$.  Then for $t\in (-\frac{1}{2},\frac{1}{2})$ the mapping $g^t(z)=z+t\varphi(z)$ is a diffeomorphism of $\ID$ to itself which extends to the identify on the boundary $\IS$.  If $f$ is a mapping of finite distortion for which $\mathsf{E}_p(f)<\infty$,  then so is $f\circ g^t$ and $f\circ g^t-f\in W^{1,1}_0(\ID)$. If $f$ is a homeomorphism,  so is $f\circ g^t$,  and so forth.  

The function $t\mapsto \mathsf{E}_p(f\circ g^t)$ is a smooth function of $t$. Thus if  $f$ is a minimiser in any reasonable class (that is we may relax the assumption that $f$ is a homeomorphism) we have
\[ \frac{d}{dt} \Big|_{t=0} \mathsf{E}_p(f\circ g^t) = 0. \]
It is a nice calculation to verify that this equation is equivalent to (\ref{1.5}).  It is interesting to note that (\ref{1.5}) implies that $\mu_f$ is constant on any open set that $|\mu_f|$ is constant.

\medskip

In this article we prove the following:
\begin{theorem}\label{thm1.1}
Let $f:\ID\to\IC$ be a finite distortion mapping that satisfies the distributional equation  (\ref{1.5}). Assume that
\begin{equation}\label{1.6}
\mathbb{K}_f=\frac{1+|\mu_f|^2}{1-|\mu_f|^2}\in L_{loc}^r(\ID),\mbox{ for some }r>p+1.
\end{equation}
Then $f$ is a local diffeomorphism from $\ID$ to $f(\ID)$.
\end{theorem}
Note that  {\em a priori} we have  $\mathbb{K}_f\in L^{p}(\ID)$, and so this result assumes slightly more than that.  We will see that the value $r$ need only be locally uniform in $\ID$. Theorem \ref{thm1.1} is essentially best possible,  see Theorem \ref{thm 1.2} below.

\medskip

The following corollary follows almost immediately:
\begin{corollary}
Let $p>1$ and $f$ be a locally quasiconformal minimiser continuous in $\overline{\ID}$ and $f|\IS=f_0$, for homeomorphic boundary values $f_0:\IS\to\IS$. Then $f:\overline{\ID}\to\overline{\ID}$ is a homeomorphism and a $C^\infty$-smooth diffeomorphism on $\ID$.
\end{corollary}
The case $p=1$ was already known,  \cite{AIMO}.  Here is an outline of the paper.\\

\S 2 will be devoted to the proof of Theorem \ref{thm1.1}.\\

\S 3  gives the following counterexample to justify the assumption of the integrability of the distortion.  A pseudo-inverse of a mapping of finite distortion is a continuous $W^{1,2}(\overline{\ID})$ surjection onto $\overline{\ID}$ with $(h\circ f)(z)=z$ almost everywhere.

\begin{theorem}\label{thm 1.2}
For each $p>1$ there is a Sobolev mapping $f:\ID\to\ID$ of finite distortion with the following properties. 
\begin{enumerate} 
\item $f\in W_{loc}^{1,\frac{2p}{p+1}}(\ID)$ and $\mathsf{E}_p(f)<\infty$,
\item the Beltrami coefficient $\mu_f$ satisfies the distributional equation (\ref{1.5}),  
\item $\IK_f \in L^p_{loc}(\ID)\setminus\bigcup_{q>p} L^q_{loc}(\ID)$, 
\item $|\ID\setminus f(\ID)|=0$ and
\item $f$ has a pseudo-inverse $h:\overline{\ID}\to \overline{\ID}$ which is monotone.
\end{enumerate}

In particular,  this mapping $f$ has $\mathsf{E}_p(f)<\infty$ and solves the distributional equation, but it cannot even be locally quasiconformal.
\end{theorem}

This mapping $f$ cannot refute Conjecture \ref{conjecture} as there are two issues.  First  it may not be a minimiser,  although it is a stationary point for smooth variations.  One therefore expects there is some as yet unexplained reason why a minimiser has higher regularity.  This situation is not uncommon though.  \\

The second issue is that concerning the boundary values of $f$.  Sometimes $f$ can be promoted to a homeomorphism of the boundary.  Briefly; by definition the pseudo-inverse $h$ has continuous boundary values $h_0:\IS\to\IS$ whose argument is continuously increasing from $[0,2\pi]$ onto $[0,2\pi]$.  As a monotone mapping $h_0$ has a countable collection of closed disjoint compression arcs $\{\alpha_i\}_{i=1}^{\infty}$ (on which it is constant).  A conformal resolution of $h_0$ is a conformal mapping $\varphi:\ID\hookrightarrow\ID$ continuous on $\overline{\ID}$ and with $\varphi(\alpha_i) \cap \IS$ a singleton,  $\varphi(\ID)$ of full measure and such that $h_0\circ(\varphi^{-1}|\IS):\IS\to\IS$ is a homeomorphism (where we understand for $\zeta\in \IS$ that $\varphi^{-1}(\zeta)$ is a set on which $h_0$ is constant). For example the Riemann mapping $\psi_\lambda:\ID \to \ID\setminus [\lambda,1)$, $0<\lambda<1$, $\psi_\lambda(0)=0$, $\psi_\lambda'(0)>0$, is a conformal resolution of the continuous map $e^{i\theta} \mapsto e^{i\eta(\theta)}$,  $\eta:[\theta_\lambda,2\pi-\theta_\lambda]\to[0,2\pi]$ a homeomorphism and  $\eta([-\theta_\lambda\leq\theta\leq\theta_\lambda])=0$.  Thus compressing the arc $\{e^{i\theta}:-\theta_\lambda\leq\theta\leq\theta_\lambda\}$ -- here $\theta_\lambda$ can be found as an explicit function of $\lambda$.   Given a conformal resolution of $h_0$ above,  the map $h\circ \varphi^{-1}$ will now be a homeomorphism on the boundary.  Notice that 
\[ \int_\ID \|D(h\circ \varphi^{-1})\|^2 = \int_{\varphi(\ID)} \|D(h\circ \varphi^{-1})\|^2 = \int_{\ID} \|Dh\|^2, \]
and that $\mu_{(h\circ \varphi^{-1})^{-1}}=\mu_f$ almost everywhere.  Thus $(h\circ \varphi^{-1})^{-1}=\varphi \circ f $ could be found in Theorem \ref{thm 1.2} with homeomorphic boundary values should $h_0$ admit a conformal resolution -- even though $\varphi \circ f $ may not be a homeomorphism. Iterating maps like $\psi_\lambda$ and some elementary normal family arguments and other considerations yields conformal mappings ``resolving'' any disjoint collection of closed arcs.  However the continuity of the composition $h_0\circ(\varphi^{-1}|\IS)$ is a real issue,  and is perhaps unlikely in general.  For certain families of compression arcs,  such as finite families,  or if $\IS\setminus \bigcup_{i} \alpha_i = \bigcup \beta_j$ is  a countable collection of open arcs, we can construct resolutions and thus in some circumstances we can promote the map of Theorem \ref{thm 1.2} to a homeomorphism.  We ask whether the integrability of the distortion of $f$ might deliver a good set of compression arcs for $h_0$? \\

\S 4 gives various equivalent conditions to imply Theorem \ref{thm1.1}. For instance if $f\in W^{1,s}_{loc}(\ID)$, $s>2$,  is a minimiser,  then the conclusions of Theorem \ref{thm1.1} are valid. See for instance Theorem \ref{thm 4.1}.\\

\S 5 we discuss the existence part of Conjecture 1.1. Theorem \ref{thm 1.2} suggests the optimal classes where one might find a minimiser. If $\{f_j\}_{j=1}^{\infty}$ is a minimising sequence of homeomorphic mappings of finite distortion, then H\"older's inequality gives the {\em a priori} bounds (see \cite{IMO})
\begin{eqnarray} \nonumber
\Big[\int_\ID \|Df_j(z)\|^{\frac{2p}{p+1}}dz\Big]^{p+1} & \leq & \int_{\ID }\IK^p(z,f_j)\; dz\; \cdot\; \Big[\int_\ID J(z,f_j)dz\Big]^p \\ & \leq & \pi^p\int_{\ID }\IK^p(z,f_j)\; dz. \label{1.7}
\end{eqnarray}
Thus there is a subsequence $f_j\rightharpoonup f$ weakly in $W^{1,\frac{2p}{p+1}}(\ID)$. Meanwhile, the sequence of inverses,  $h_j=f_j^{-1}$, satisfies
\begin{align} 
\int_\ID\|Dh_j(w)\|^2dw&=\int_\ID\IK(w,h_j)J(w,h_j)dw\notag\\
&\leq\int_{\ID }\IK^p(w,h_j)J(w,h_j)dw=\int_{\ID }\IK^p(z,f_j)dz,\label{1.8}
\end{align}
so they converge in $W^{1,2}(\ID)$. Note the change of variables formula follows from \cite{H,KM}. Such a sequence $h_j$ converges to a continuous function $h$ locally uniformly in $\ID$ \cite{GV,IM}.  In fact in \cite{IKO} it is proved that $h$ will be locally Lipschitz. However, on the $f$ side, functions in $W^{1,\frac{2p}{p+1}}(\ID)$ are not usually continuous.  The continuity of $f$ is the main obstacle to it being a homeomorphism.
To address this  problem we define a larger space $\mathsf{H}_p$ in which the minimising sequence converges to an inner-variational minimiser, thus $f$ will satisfy equation (\ref{1.5}). Furthermore, if it satisfies the hypotheses of Theorem \ref{1.1}, then it is a diffeomorphism onto its image, and a minimiser in $\mathsf{F}_p$ if it is continuous in $\overline{\ID}$.\\

\S 6  gives an analytic condition which implies the minimiser will be a homeomorphism. Roughly $f\in W^{1,2}(\ID)$ is enough.\\

\S 7  discusses the limiting regimes $p\to\infty$ and $p\to 1$.  The case $p=1$ is by now well-known \cite{AIMO}.  But our direct methods here fail for it as we do not get a uniform elliptic estimate and we resort to an alternative approach.  We prove that as $p\to 1$ the psuedo-inverses of minimisers converge locally uniformly to the harmonic mapping.  When $p\to\infty$ we show the local uniform limit exists and is an extremal quasiconformal mapping and identify when the approximating sequence is a Hamilton sequence,  making this limit a uniquely extremal Teichm\"uller mapping for its boundary values.  \\

\S 8 for each $p\geq 1$ we give examples with non-constant Ahlfors-Hopf differential $\phi=\alpha \,z^{-2}$.  These are them used to show that there are quasisymmetric mappings $f_0:\IS\to\IS$ with diffeomorphic extensions $f:\ID\to\ID$ which are minimisers of $\mathsf{E}_p$ and are uniquely so even in the larger class of mappings of finite distortion with boundary values $f_0$.

\section{Diffeomorphisms ; proof of Theorem \ref{1.1}.}
We rewrite the distributional equation as
\[
2p\int_\ID \IK_f^p\frac{\overline{\mu_f}}{1+|\mu_f|^2}\varphi_\zbar dz=\int_\ID \big(\IK_f^p-1\big)\varphi_zdz,\quad\forall\varphi\in C_0^\infty(\ID ).
\]
Let $r>p+1$ and $s=r/p>1+\frac{1}{p}$.  Note that the following argument is entirely local.
\begin{lemma}\label{2.1}
There is an $F\in W_{loc}^{1,s}(\ID )$ such that
\begin{equation}
F_z=2p\IK_f^p\frac{\overline{\mu_f}}{1+|\mu_f|^2},\quad F_\zbar =\IK_f^p-1.
\end{equation}
\end{lemma}
\noindent{\bf Proof.}
We write $a(z)=2p\IK_f^p\frac{\overline{\mu_f}}{1+|\mu_f|^2}$, $b(z)=\IK_f^p-1$. By assumption they are both in $L_{loc}^s(\ID)$. Let $0<r<1$. Define
\[
a^r(z)=\begin{cases}
a(z),&z\in D_r,\\
0,&z\in\IC\setminus D_r,
\end{cases} \quad {\rm and} \quad
b^r(z)=\begin{cases}
b(z),&z\in D_r,\\
0,&z\in\IC\setminus D_r.
\end{cases}.
\]
Then $a^r,b^r\in L^s(\IC)$. Define $G=\mathcal{C}(b^r)$, $H=\mathcal{C}^*(a^r)$, where $\mathcal{C}$ is the Cauchy transform, and $\mathcal{C}^*$ is its conjugate defined by $\mathcal{C}^*\eta=\overline{\mathcal{C}\overline{\eta}}$. Note that both $G$ and $H$ are in $W^{1,s}(\IC)$. Then
\[
\int_{D_r}G\varphi_{z\zbar }=-\int_{D_r}G_\zbar \varphi_z=-\int_{D_r}H_z\varphi_\zbar =\int_{D_r}H\varphi_{z\zbar }.
\]
It follows from Weyl's lemma that in $D_r$, $\phi:=G-H$ is harmonic, and then $\phi_z=G_z-H_z$ is holomorphic. Let $\Phi$ be an anti-derivative of $\phi_z$, and set
$F^r=G-\Phi$. Then in $D_r$ we have
\[
F^r_z=G_z-\phi_z=H_z=a,\quad F^r_\zbar =G_\zbar=b.
\]
We consider $r<R<1$. By the same process we get an $F^R$ defined in $D_R$, and in $D_r$ we have $F^R_z=F^r_z=a$, $F^R_\zbar=F^r_\zbar=b$, thus $F^r-F^R=c_0$, a constant. If we redefine $F^R$ by $F^R+c_0$, then we have $F^R=F^r$ in $D_r$. Now the function $F(z):=F^{(1+r)/2}(z)$ for $|z|<r$ is well defined in $\ID$ and it also satisfies the conditions required.
\hfill $\Box$

\bigskip

Now for $F$ as in Lemma \ref{2.1}, we write $t=|\mu_f|$ and calculate that
\[
|F_z|=2p\Big(\frac{1+t^2}{1-t^2}\Big)^p\frac{t}{1+t^2},\quad F_\zbar =\Big(\frac{1+t^2}{1-t^2}\Big)^p-1.
\]
Then
\[
t=\sqrt{\frac{(F_\zbar +1)^\frac{1}{p}-1}{(F_\zbar +1)^\frac{1}{p}+1}}.
\]
We can therefore write
\begin{equation}\label{2.2}
|F_z|=a_p(F_\zbar ),
\end{equation}
where
\begin{equation}\label{2.3}
a_p(s)=p(s+1)^\frac{p-1}{p}\sqrt{(s+1)^\frac{2}{p}-1}.
\end{equation}
We compute
\begin{equation}\label{2.4}
a_p^\prime(s)=\frac{p(s+1)^\frac{1}{p}-(p-1)(s+1)^{-\frac{1}{p}}}{\sqrt{(s+1)^\frac{2}{p}-1}},
\end{equation}
and find
\[
a_p^\prime(+\infty)=p.
\]
Now
\begin{eqnarray*}
a_p(s)-ps&= &p(s+1)^\frac{p-1}{p}\sqrt{(s+1)^\frac{2}{p}-1}-ps =p\Big(\sqrt{(s+1)^2-(s+1)^\frac{2(p-1)}{p}}-s\Big)\\
&=& p\frac{2s+1-(s+1)^\frac{2(p-1)}{p}}{\sqrt{(s+1)^2-(s+1)^\frac{2(p-1)}{p}}+s} = \begin{cases}
\mathcal{O}(1),&1<p\leq2,\\
\mathcal{O}(s^{1-\frac{2}{p}}),&p>2,
\end{cases}
\end{eqnarray*}
as $s\to\infty$. So we conclude
\[
a_p(s)=ps+\mathcal{O}(s^\alpha),\quad\alpha=
\begin{cases}
0,&1<p\leq2;\\
1-\frac{2}{p},&p>2.
\end{cases}
\]
Precisely,
$ |F_z|=pF_\zbar +\mathcal{O}(F_\zbar ^\alpha)$.
Set
\[
\kappa=\frac{1}{p}\frac{\overline{F_z}}{|F_z|}.
\]
Then
\begin{equation}\label{2.5}
F_\zbar =\kappa F_z+\mathcal{O}(F_\zbar ^\alpha),\quad|\kappa|=\frac{1}{p},\quad\alpha=
\begin{cases}
0,&1<p\leq2,\\
1-\frac{2}{p},&p>2.
\end{cases}
\end{equation}
Note that $\alpha<1$.
\begin{lemma}
Let $F\in W_{loc}^{1,s}(\ID )$, $s>1+\frac{1}{p}$ satisfy equation (\ref{2.5}). Then $F\in W_{loc}^{1,2}(\ID )$.
\end{lemma}
\noindent{\bf Proof.} Let $\eta\in C_0^\infty(\ID )$. Then,
\begin{eqnarray*}
(\eta F)_\zbar &=&\eta_\zbar F+\eta F_\zbar = \eta_\zbar F+\eta\kappa F_z+\mathcal{O}(\eta F_\zbar ^\alpha)\\
&=& \eta_\zbar F+\kappa[(\eta F)_z-\eta_zF]+\mathcal{O}(\eta F_\zbar ^\alpha).
\end{eqnarray*}
Note here we have $F\in L_{loc}^{s^*}(\ID )$, where $s^*>2$ is the Sobolev conjugate of $s$. For $r=\min\{s^*,\frac{s}{\alpha}\}$,
\begin{equation}\label{2.6}
(\mathbf{I}-\kappa\mathcal{S})(\eta F)_\zbar =\eta_\zbar F-\kappa\eta_zF+\mathcal{O}(F_\zbar ^\alpha)\in L^r(\IC ).
\end{equation}
Meanwhile, $s>1+\frac{1}{p}$ implies that $\mathbf{I}-\kappa\mathcal{S}$ is invertible on $\eta F$, see \cite[Theorem 14.0.4]{AIM}. Thus we have $(\eta F)_\zbar \in L^r(\IC )$. If $\frac{s}{\alpha}\geq2$, then the claim follows; otherwise we have $F\in W_{loc}^{1,\frac{s}{\alpha}}(\ID )$, which puts the right-hand side of (\ref{2.6}) in $L^{r^\prime}(\IC )$ for $r^\prime=\min\{(\frac{s}{\alpha})^*,\frac{s}{\alpha^2}\}$. Again, if $\frac{s}{\alpha^2}\geq2$, then the proof is complete; otherwise we keep iterating until $\frac{s}{\alpha^n}\geq2$, and this completes the proof.\hfill $\Box$

\bigskip

\begin{lemma}
Suppose $F\in W^{1,2}_{loc}(\ID)$ satisfies (\ref{2.2})-(\ref{2.3}). Then $F$ is smooth.
\end{lemma}
\noindent{\bf Proof.} We start from (\ref{2.4}):
\[
a_p^\prime(s)=\frac{p(s+1)^\frac{1}{p}-(p-1)(s+1)^{-\frac{1}{p}}}{\sqrt{(s+1)^\frac{2}{p}-1}}.
\]
Write $
x=(s+1)^\frac{1}{p}$, and 
\[
c_p(x)=\frac{px-(p-1)x^{-1}}{\sqrt{x^2-1}}, \;\;\;\; 
c_p^\prime(x)=\frac{(p-2)x^2-(p-1)}{x^2(x^2-1)^\frac{3}{2}}.
\]
Thus for $1<p\leq2$, $c_p(x)$ has a minimum value $c_p(+\infty)=p$,  while for $p>2$, $c_p(x)$ attains its minimum
$
c_p(\sqrt{(p-1)/(p-2)}\,)=2\sqrt{p-1}.
$
Hence for any $p>1$,
 we have
\[ a_p^\prime(s)\geq M_p = \left\{\begin{array}{ll}  p, &1<p\leq 2. \\
2\sqrt{p-1}, & p\geq 2. \end{array} \right.   \]
Since $M_p>1$, $a_p$ is increasing and we can now write (\ref{2.2}) as follows.
\begin{equation}\label{2.7}
F_\zbar =\mathcal{A}_p(|F_z|), \hskip10pt |\mathcal{A}_p^\prime|\leq k_p = \frac{1}{M_p} <1.
\end{equation}
This is an elliptic equation in the sense of [2, Definition 7.7.1] and then it follows from the Caccioppoli-type estimates ([2, Theorem 5.4.2], also see [2, Theorem 8.7.1]) that $F\in W^{2,2}_{loc}(\ID )$. We next consider the function
\[
|F_z|^2=a_p^2(F_\zbar )=b_p(F_\zbar ),
\]
where
\begin{eqnarray*}
b_p(s)&=&p^2[(s+1)^2-(s+1)^\frac{2(p-1)}{p}],\\
b_p^\prime(s)&=&p^2[2(s+1)-\frac{2(p-1)}{p}(s+1)^\frac{p-2}{p}],\\
b_p^{\prime\prime}(s)&=&p^2[2-\frac{2(p-1)(p-2)}{p^2}(s+1)^{-\frac{2}{p}}]\geq0.
\end{eqnarray*}
Then  $
\min_{s\geq0}b_p^\prime(s)=b_p^\prime(0)=2p>1
$,  thus $b$ is invertible and we can write $\mathcal{B}_p=b_p^{-1}$, and
\begin{equation}\label{2.8}
F_\zbar =\mathcal{B}_p(|F_z|^2),
\end{equation}
Note here $\mathcal{B}_p(t^2)=\mathcal{A}_p(t)$, thus
\begin{equation}
\mathcal{A}_p^\prime(t)=2t\mathcal{B}_p^\prime(t^2).
\end{equation}
As $F\in W^{2,2}_{loc}(\ID )$, we may differentiate both sides of (\ref{2.7}) by $x$, and get
\begin{equation}\label{2.9}
(F_x)_\zbar =\mathcal{B}_p^\prime(|F_z|^2)\overline{F_z}(F_x)_z+\mathcal{B}_p^\prime(|F_z|^2)F_z\overline{(F_x)_z},
\end{equation}
where it follows from (\ref{2.8}) that
\[
\mathcal{B}_p^\prime(|F_z|^2)|\overline{F_z}|+|\mathcal{B}_p^\prime(|F_z|^2)|F_z|\leq\mathcal{A}_p^\prime(|F_z|)\leq k_p.
\]
Thus (\ref{2.9}) is again an elliptic equation for the function $F_x$, thus $F_x\in W^{2,2}_{loc}(\ID )$. The same argument applies for the function $F_y$, so that $F\in W_{loc}^{3,2}(\ID )$. So we can differentiate (\ref{2.9}) again,  
\[
(F_{xx})_\zbar =\mathcal{B}_p^\prime(|F_z|^2)\overline{F_z}(F_{xx})_z+\mathcal{B}_p^\prime(|F_z|^2)F_z\overline{(F_{xx})_z}+\Phi(z),
\]
where $\Phi(z)$ is composed of lower-order terms, so the equation is again elliptic. Now the argument is inductive and so we conclude that $F$ is smooth.\hfill $\Box$

\bigskip

We now observe the smoothness of $\mu_f$ follows. In fact by (\ref{2.1}),
\begin{equation}
\mu_f=\frac{\overline{F_z}}{|F_z|}\sqrt{\frac{(F_\zbar +1)^\frac{1}{p}-1}{(F_\zbar +1)^\frac{1}{p}+1}}\\
=\frac{\overline{F_z}}{p[(F_\zbar+1)+(F_\zbar+1)^\frac{p-1}{p}]},
\end{equation}
which is smooth as $F_\zbar\geq 0$.  We now need the following lemma.

\begin{lemma}
If a finite distortion function $f:\Omega\to\IC$ has smooth Beltrami coefficient $\mu$ and $p$-integrable distortion $\IK(z,f)$,  
\[
\int_\Omega\IK^p(z,f)dz<\infty,
\]
for some $p\geq1$, then $|\mu_f|<1$ in $\Omega$.
\end{lemma}
\noindent{\bf Proof.} Suppose $|\mu_f(z_0)|=1$, for some $z_0\in\Omega$. For notational ease we set $z_0=0$ and consider the function $|\mu_f|$ to be smooth in a disk $D(0,\delta)$. As $|\mu_f|\leq1$, we have $|\mu_f|_x(0)=|\mu_f|_y(0)=0$. From Taylor's expansion,
\[
|\mu_f(z)|\geq1-M|z|^2,\quad z\in D(0,\delta),
\]
where $M=\sup_{z\in D(0,\delta)}|\,\nabla^2|\mu_f|\,|<\infty.$
Then,
\[ 
\int_\Omega\Big(\frac{1+|\mu_f|^2}{1-|\mu_f|^2}\Big)^p \geq\frac{1}{2^p}\int_\Omega\frac{1}{(1-|\mu_f|)^p} \geq\frac{1}{(2M)^p}\int_{D(0,\delta)}\frac{1}{|z|^{2p}}=\infty,
\]
which gives the contradiction.\hfill $\Box$

\bigskip

\noindent{\bf Proof of Theorem 1.1.} Let $\Omega\subset\subset\ID $ be compactly contained. By Lemma 2.4, there is a $k$ such that
\[
|\mu_f(z)|\leq k<1,\quad\forall z\in\Omega.
\]
So $f$ is locally quasiregular in $\ID$, with a smooth Beltrami coefficient $\mu_f$. Such a function is locally diffeomorphic.\hfill $\Box$

\section{A counterexample; proof of Theorem \ref{thm 1.2}.}

We begin with the following lemma.  We refer to \cite[\S 5.4]{AIM} for discussion of the critical interval.
\begin{lemma}[Existence of Solution]
Let $\Omega\subset\IC$ be a planar domain of finite measure and $\mathcal{H}:\Omega\times \IC$ be measurable in the first coordinate and $H(z,0)\equiv 0$, for every $z\in\Omega$. Suppose that there is a $k<1$ and $\alpha<1$ such that for all $z\in\Omega$ and $\zeta\in\IC$,
\[
\mathcal{H}(z,\xi)\leq k|\xi|+C |\xi|^\alpha+|h(z)|,
\]
where $C>1$ is a constant, and $h\in L^s(\Omega)$ for some $s\in(Q(k),P(k))$, where $(Q(k),P(k))$ is the critical interval for $q$ such that the operator norm $\mathbf{S}_q$ of $\mathcal{S}:L^q(\IC)\rightarrow L^q(\IC)$ satisfies
\[
k\mathbf{S}_q<1.
\] 
Then, for every $q\in(Q(k),s]$, there is an $f\in W^{1,q}(\Omega)$ that satisfies
\begin{equation}\label{3.1}
f_\zbar=\mathcal{H}(z,f_z), \hskip20pt a.e. \;\; z\in\Omega.
\end{equation}
\end{lemma}
\noindent{\bf Proof.} We extend $\mathcal{H}$ by zero outside $\Omega$,  $\mathcal{H}(z,\zeta)=0$,  in  $ \IC\setminus \Omega$ and work in $\IC$.  
Let $g\in L^q(\Omega)$ and extend $g$ by $0$ to an element of $L^p(\IC)$. Then H\"older's inequality gives
\begin{eqnarray*}
\|  |\mathcal{S}g|^{\alpha} \|_{L^q(\Omega)} & = & \Big( \int_\Omega |\mathcal{S}g|^{\alpha q}\Big)^{1/q} \leq |\Omega|^{\frac{1}{q(1-\alpha)}} \Big(\int_\Omega |\mathcal{S}g|^q\Big)^{\alpha/q}  \\
&\leq & |\Omega|^{\frac{1}{q(1-\alpha)}} \Big(\int_\IC|\mathcal{S}g|^q\Big)^{\alpha/q} = |\Omega|^{\frac{1}{q(1-\alpha)}} \| \mathcal{S}g\|_{L^q(\IC)}^{\alpha} \\
& \leq &  |\Omega|^{\frac{1}{q(1-\alpha)}} \mathbf{S}_q^\alpha \|g\|_{L^q(\IC)}^{\alpha} =  |\Omega|^{\frac{1}{q(1-\alpha)}} \mathbf{S}_q^\alpha \|g\|_{L^q(\Omega)}^{\alpha} .
\end{eqnarray*}
We now compute as follows.
\begin{eqnarray*}
\| \mathcal{H}(z,\mathcal{S}g)\|_{L^q(\Omega)} & \leq  &k \|\mathcal{S}g\|_{L^q(\Omega)} +C \| |\mathcal{S}g|^\alpha \|_{L^q(\Omega)}+\| h \|_{L^q(\Omega)}\\
& \leq  &k\mathbf{S}_q \|g\|_{L^q(\IC)}  +C |\Omega|^{\frac{1}{q(1-\alpha)}} \mathbf{S}_q^\alpha \|g\|_{L^q(\Omega)}^{\alpha}+\| h \|_{L^q(\Omega)} \\
& =  &k\mathbf{S}_q \|g\|_{L^q(\Omega)}  +C |\Omega|^{\frac{1}{q(1-\alpha)}} \mathbf{S}_q^\alpha \|g\|_{L^q(\Omega)}^{\alpha} +\| h \|_{L^q(\Omega)}.
\end{eqnarray*}
Hence if $\|g\|_{L^q(\Omega)}\leq R$ we have
\begin{eqnarray*}
\| \mathcal{H}(z,\mathcal{S}g)\|_{L^q(\Omega)} & \leq  & k\mathbf{S}_q R  +C |\Omega|^{\frac{1}{q(1-\alpha)}} \mathbf{S}_q^\alpha R^{\alpha} +\| h \|_{L^q(\Omega)}.
\end{eqnarray*}
Thus as soon as $k\mathbf{S}_q<1$ we can find a sufficiently large $R$ so that 
\[ \|\mathcal{H}(z,\mathcal{S}g)\|_{L^q(\Omega)} < R \] and we are in a position to apply  the Schauder fixed-point theorem to find a $g\in L^q(\Omega)$ such that
$g=\mathcal{H}(z,\mathcal{S}g)$.
Then the mapping $f=\mathcal{C}g$ is a $W^{1,q}(\IC)$ solution to (\ref{3.1}).\hfill $\Box$

\bigskip

\noindent In \cite{Astala},  K. Astala made us aware of a counterexample to a potential generalisation of the super-regularity theorem for autonomous Beltrami systems, \cite{Martin}.
\begin{lemma}\label{lem 3.1}
There is a function $G\in W^{1,1}_{loc}(\ID)$,  but not in $W^{1,q}_{loc}(\ID)$ for any $q>1$, that satisfies
\[
G_\zbar(z)=\frac{1}{p}|G_z(z)|,\quad a.e.\quad z\in\ID.
\]
\end{lemma}
The function $G$ cannot be a mapping of finite distortion despite the fact that $\IK(z,f)\equiv \frac{p^2+1}{p^2-1}$ as $J(z,f)=\big(1-\frac{1}{p}\big)|f_z|^2 \in L^{1/2}_{loc}(\ID)$ and no better.

\medskip

We start with the $G$ provided by Lemma \ref{lem 3.1}. Set $a\in W_{loc}^{1,1}(\ID)$, and $F=G+a$. We wish to find an equation for $a$ which ensures $F$ satisfies (\ref{2.2})-(\ref{2.3}). After an elementary computation this can be written as
\[
F_{\bar{z}}=\mathcal{A}_p(|F_z|)=\frac{1}{p}|F_z|+\mathcal{O}(|F_\zbar|^\alpha),\quad\alpha=
\begin{cases}
0,&1<p\leq2;\\
1-\frac{2}{p},&p>2.
\end{cases}
\]
In fact,
\begin{eqnarray*}
F_\zbar&=&G_\zbar+a_\zbar=\frac{1}{p}|G_z|+a_\zbar\\
&=&\mathcal{A}_p(|F_z|)-\mathcal{A}_p(|G_z+a_z|)+\frac{1}{p}|G_z|+a_\zbar.
\end{eqnarray*}
Thus we require that almost everywhere in $\ID$,
\begin{equation}\label{3.2}
a_\zbar=\mathcal{A}_p(|G_z+a_z|)-\frac{1}{p}|G_z|
\end{equation}
We then calculate that
\begin{eqnarray*}
a_\zbar&=&\mathcal{A}_p(|G_z+a_z|)-\frac{1}{p}|G_z| = \frac{1}{p}|G_z+a_z|-\frac{1}{p}|G_z|+\mathcal{O}(|G_\zbar+a_\zbar|^\alpha) \\  
&\leq & \frac{1}{p}|a_z|+C_1|G_\zbar+a_\zbar|^\alpha+C_2 \leq   \frac{1}{p}|a_z|+C_1|a_\zbar|^\alpha+C_1|G_\zbar|^\alpha+C_2 
\end{eqnarray*}
This can be written as
\begin{equation}
g=\mathcal{H}(z,g),\hskip20pt a.e. \;\; z\in\ID,
\end{equation}
where
\begin{align*}
\mathcal{H}(z,g)\leq\frac{1}{p}|\mathcal{S}g|+C|g|^\alpha+|h(z)|,
\end{align*}
and $h\in L^\frac{1}{\alpha}(\ID)$. Here we need to check that the critical exponent
$Q(1/p)<\frac{1}{\alpha}$.
This recalls  Iwaniec's conjecture \cite{Iwaniec} that $\mathbf{S}_q=\frac{1}{q-1}$, for $1<q<2$. Although this has not been completely proved, Nazarov-Volberg \cite{NV} showed that
\[
\mathbf{S}_q=\frac{C_q}{q-1},\quad 1\leq C_q\leq 2.
\]
In fact today the current best bound known is $C_q\leq 1.575$, see Ba\~nuelos-Janakiraman \cite{BJ}. Nevertheless $C_q=2$ is enough for our purpose:
\[
\frac{1}{p}\mathbf{S}_{\frac{1}{\alpha}}\leq\frac{1}{p}\cdot\frac{2}{\frac{1}{\alpha}-1}=\frac{p-2}{p}<1.
\]
Now, by Lemma 3.1, for any $q\in(Q(\frac{1}{p}),\frac{1}{\alpha}]$, there is a $g\in L^q(\IC)$ that satisfies (3.3), thus $\mathcal{C}g$ is a $W^{1,q}(\IC)$ solution to (3.2). In particular, we can choose $a$ as a $W^{1,\frac{1}{\alpha}}(\IC)$ solution. Consider
\[
F=G+a.
\]
Then $F$ is a $W^{1,1}_{loc}(\ID)$ function but not in $W^{1,q}_{loc}(\ID)$ for any $q>1$. This establishes the following:

\begin{lemma}
There is a function $F\in W^{1,1}_{loc}(\ID)$,  but not in $W^{1,q}_{loc}(\ID)$, for any $q>1$ such that
\[
|F_z|=a_p(F_\zbar),\hskip20pt a.e. \;\; z\in\ID.
\]
where $a_p$ is defined at (\ref{2.3}).
\end{lemma}

Next, for the function $F$ of Lemma 3.3, we can set
\[
\mu=\frac{\overline{F_z}}{|F_z|}\sqrt{(\frac{F_\zbar+1)^\frac{1}{p}-1}{(F_\zbar+1)^\frac{1}{p}+1}}.
\]
This gives $
F_\zbar=\Big(\frac{1+|\mu|^2}{1-|\mu|^2}\Big)^p-1$, and
\[
|F_z|=a_p(F_\zbar)=2p\Big(\frac{1+|\mu|^2}{1-|\mu|^2}\Big)^p\frac{|\mu|}{1+|\mu|^2},
\]
so that
\[
F_z=2p\Big(\frac{1+|\mu|^2}{1-|\mu|^2}\Big)^p\frac{\mubar}{1+|\mu|^2}.
\]
Then
\[
\int_\ID F_z\varphi_\zbar=-\int_\ID F\varphi_{z\zbar}=\int_\ID F_\zbar\varphi_z,\quad\forall\varphi\in C_0^\infty(\ID).
\]
This proves that $\mu$ satisfies (\ref{1.5}), with
\[
\frac{1+|\mu|^2}{1-|\mu|^2}\in L_{loc}^p(\ID)\setminus \bigcup_{q>p}L_{loc}^q(\ID).
\]

We require the following lemma. 
\begin{lemma}[Modulus lemma]
Let $f:\ID\to\ID$ be a homeomorphism of finite distortion, $\IK(z,f)\in L^1(\ID)$, and $f(0)=0$. Then, for any disk $D_R\subset\ID$, where $R\in(0,1)$, there is an $R^\prime\in(0,1)$ such that $f(D_R)\subset\overline{D_{R^\prime}}$. Furthermore, $R^\prime$ depends only on $R$ and $\|\IK(z,f)\|_{L^1(\ID)}$.
\end{lemma}
\noindent{\bf Proof.} We consider the annulus $A(R,1)$ and its image $f(A(R,1))$. Our aim is to get an estimate of the modulus of the ring $f(A(R,1))$. Let $\gamma$ be a path that connects the two boundaries of $f(A(R,1))$, and set $h=f^{-1}$. Then
\[
\frac{1}{1-R}\int_\gamma\|Dh\| |d\gamma|\geq\frac{1}{1-R}\int_{h(\gamma)}|d(h(\gamma))|\geq1,
\]
as $h(\gamma)$ is a path connecting the two circles $|z|=R$ and $|z|=1$. Thus $\frac{1}{1-R}\|Dh\|$ is an admissible function on $f(A(R,1))$. Also, we can compute the area integral
\[
M:=\frac{1}{(1-R)^2}\int_{f(A(R,1))}\|Dh(w)\|^2dw=\frac{1}{(1-R)^2}\int_{A(R,1)}\IK(z,f)dz<\infty.
\]
So
\[
\mbox{Mod}(f(A(R,1)))=\sup_{\tilde{\rho}\in\Gamma}\frac{1}{\int_{f(A(R,1))}\tilde{\rho}^2}\geq\frac{1}{M},
\]
where $\Gamma$ is the collection of all admissible functions on $f(R,1)$. Now we have a lower bound for the Modulus of $f(A(R,1))$. Also we have the assumption $f(0)=0$, which implies
\[
\min_{|z|=R}|f(z)-1|>0.
\]
So $R^\prime:=1-\min_{|z|=R}|f(z)-1|$ satisfies the requirements.\hfill $\Box$

\bigskip

\begin{theorem}
Let $\mu:\ID\to\overline{\ID}$ be measurable, and
\[
\int_\ID\Big(\frac{1+|\mu|^2}{1-|\mu|^2}\Big)^p<\infty,\quad p\geq1.
\]
Then there exists a finite distortion function $f\in W^{1,\frac{2p}{p+1}}(\ID)$ that satisfies the following conditions:

\begin{itemize} 
\item $f$ satisfies the Beltrami coefficient
\[
f_z=\mu f_\zbar
\]
almost everywhere in $\ID$.
\item There is a finite distortion function $h\in C(\overline{\ID})\cap W^{1,2}(\ID)$, $h$ is monotone in $\overline{\ID}$, and
\[
\int_\ID\mathbb{K}^p(w,h)J(w,h)=\int_\ID\Big(\frac{1+|\mu|^2}{1-|\mu|^2}\Big)^p.
\] 
\item There is a measurable set $X\subset\ID$ such that $|\ID-X|=0$, $h\circ f(z)=z$ for every $z\in X$, and $J(w,h)=0$ for almost every $w\in\ID-f(X)$.
\end{itemize}
\end{theorem}
\noindent{\bf Proof.} Set
\[
\mu^m(z)=\begin{cases}
\mu(z),&\mbox{ if }|\mu(z)|\leq1-\frac{1}{m};\\
(1-\frac{1}{m})\frac{\mu(z)}{|\mu(z)|},&\mbox{otherwise}.
\end{cases}
\]
For each $m$ there is a quasiconformal mapping $f^m:\overline{\ID}\to\overline{\ID}$ such that $\mu_{f^m}=\mu^m$ almost everywhere  in $\ID$. Also note that
\[
\Big[\int_\ID\|Df^m\|^\frac{2p}{p+1}\Big]^{p+1}\leq\pi^p\int_\ID\Big(\frac{1+|\mu^m|^2}{1-|\mu^m|^2}\Big)^p\leq\pi^p\int_\ID\Big(\frac{1+|\mu|^2}{1-|\mu|^2}\Big)^p.
\]
Up to a subsequence there is a limit function $f$ such that $f^m\rightharpoonup f$ in $W^{1,\frac{2p}{p+1}}(\ID)$. We next show that $\mu_f=\mu$. In fact, $\mu^m\to\mu$ pointwise. Let $\phi\in C_0^\infty(\ID)$. Then
\begin{eqnarray*}
\lefteqn{\Big|\int_\ID\phi(\mu f_z-\mu^mf_z^m)\Big|}\\&\leq &\Big|\int_\ID\phi(\mu f_z-\mu f_z^m)\Big|+\Big|\int_\ID\phi(\mu f_z^m-\mu^mf_z^m)\Big|\\
&\leq& \Big|\int_\ID\phi\mu(f_z-f_z^m)\Big|+\|\phi\|_{L^\infty(\ID)}\|f^m_z\|_{L^q(\ID)}\|\mu-\mu^m\|_{L^{q^\ast}(\ID)} 
 \to0,
\end{eqnarray*}
where $q^\ast$ is the H\"older conjugate of $q=\frac{2p}{p+1}$. So
\begin{align*}
\int_\ID\phi(f_\zbar-\mu f_z)=\lim_{m\to\infty}\int_\ID\phi(f^m_\zbar-\mu^mf^m_z)=0.
\end{align*}
This proves $f_\zbar=\mu f_z$ almost everywhere in $\ID$.\\
For the rest claims we consider the inverse sequence $h^m=(f^m)^{-1}$. By the modulus lemma, we can choose any $0<r<1$, then $|h^m(w)|\geq r'>0$ for $w\in\overline{\ID}-D_r$. Thus we can extend $h_m$ to $D(0,\frac{1}{r})$ by defining 
\[ h^m(w)=\frac{1}{\overline{h^m(\frac{1}{\wbar})}} ,\] if $w\in D(0,\frac{1}{r})-\ID$. Now $h^m$ is a sequence of finite distortion homeomorphisms with uniformly bounded $\|Dh^m\|_{L^2(D(0,\frac{1}{r}))}$, so it converges uniformly to a monotone limit $h$ in $\overline{\ID}$, and the other requisite properties are preserved under this convergence. The exact details here are discussed fully in \S 5  below. \hfill $\Box$

\section{Alternative conditions implying smoothness.}
In \cite{IMO} it is also proved that the inverse $h=f^{-1}$ of a homeomorphic minimiser $f$ of the $L^p$ problem satisfies the following inner-variational equation
\begin{equation}\label{4.1}
\int_\ID \IK^{p-1}(w,h)h_w\overline{h_\wbar }\varphi_\wbar dw=0,\quad\forall\varphi\in C_0^\infty(\ID ).
\end{equation}
In fact this is also true of diffeomorphic critical points of $\mathsf{E}_p$ as well.  We remark that the kernel in (\ref{4.1}) can be written as
\begin{equation}\label{Alfors}
\Phi(w)=\IK^{p-1}(w,h)h_w\overline{h_\wbar }=\IK^p(w,h)J(w,h)\frac{\overline{\mu_h}}{1+|\mu_h|^2}\in L^1(\ID ).
\end{equation}
By Weyl's lemma $\Phi$ is a holomorphic function. In particular, $\Phi\in L_{loc}^\infty(\ID)$. Ahlfors first realised this in \cite[\S 4, pg. 45]{Ahlfors} and so we call the function $\Phi$ defined at (\ref{Alfors}) the Ahlfors-Hopf differential.

\medskip

To get smoothness of a minimiser $f$, our approach requires (\ref{1.4}):
\[
\IK(z,f)\in L^r_{loc}(\ID),\quad r>p+1.
\]
We seek alternatives conditions for this, with the help of the Ahlfor's-Hopf differential. First
\begin{eqnarray*}
\Phi(f)&=&\IK^{p-1}(z,f)\frac{\overline{f_zf_\zbar}}{J_f^2} =
 \IK^p(z,f)\frac{\overline{\mu_ff_z}}{J_ff_z(1+|\mu_f|^2)}\\
&=&\IK^{p+1}(z,f)\frac{\overline{\mu_f}}{(1+|\mu_f|^2)^2}\frac{1}{f_z^2}.
\end{eqnarray*}
As $\Phi(f)\in L_{loc}^\infty(\ID)$, we observe
\[
\IK^p=\mathcal{O}(J_f),\quad\IK^{p+1}=\mathcal{O}(|f_z^2|)
\]
So (\ref{1.5}) is satisfied if either $f\in W^{1,r}_{loc}(\ID)$ for some $r>2$ or $J(z,f)\in L^s_{loc}(\ID)$ for some $s>1+\frac{1}{p}$. On the other hand we also have
\begin{align*}
\int_\ID\IK^r(z,f)dz&=\int_\ID\IK^r(w,h)J(w,h)dw\\
&=\int_\ID\IK^{r-p}(w,h)\Big(\IK^p(w,h)J(w,h)\Big)dw,
\end{align*}
or
\begin{align*}
\int_\ID\IK^r(z,f)dz&=\int_\ID\IK^r(w,h)J(w,h)dw\\
&=\int_\ID\Big(\IK^p(w,h)J(w,h)\Big)^\frac{r}{p}J^{1-\frac{r}{p}}(w,h)dw.
\end{align*}
Then, to get $\IK_f\in L_{loc}^r(\ID)$, $r>p+1$, we need $\IK_h^rJ_h\in L_{loc}^1(\ID)$, for $r>p+1$; or $\IK_h^s\in L_{loc}^1(\ID)$, for $s>1$; or $J_h^{-\varepsilon}\in L_{loc}^1(\ID)$, for $\varepsilon>\frac{1}{p}$. We collect all of the conditions and record as follows:\\

\begin{theorem}\label{thm 4.1}
Let $f$ be a minimiser in $\mathsf{F}_p$, and $h=f^{-1}$. Then $f$ is a diffeomorphism from $\overline{\ID}$ to $\overline{\ID}$, if any one of the following conditions is satisfied:\\
\begin{enumerate}
\item
$\IK(z,f)\in L^r_{loc}(\ID),\quad r>p+1$;
\item
$
f\in W^{1,s}_{loc}(\ID),\quad s>2;
$
\item 
$J^{1+\varepsilon}(z,f)\in L^1_{loc}(\ID),\quad\varepsilon>\frac{1}{p}$;
\item
$\IK^r(w,h)J(w,h)\in L^1_{loc}(\ID),\quad r>p+1$;
\item
$\IK(w,h)\in L^s_{loc}(\ID),\quad s>1$;
\item
$J^{-\varepsilon}(w,h)\in L^1_{loc}(\ID),\quad\varepsilon>\frac{1}{p}$.
\end{enumerate}
\end{theorem}

\section{The enlarged space}
As we discussed in the introduction, the space $\mathsf{F}_p$ is not closed under weak convergence, so we do not necessarily have  existence of a minimiser in this space for which we might seek improved regularity. 

In this section we will enlarge the space $\mathsf{F}_p$ so as to be closed and consequently identify an inner variational minimiser in this larger space.  This minimiser satisfies the variational equation (\ref{1.5}). In particular, if any of the conditions in Theorem \ref{4.1} holds, then it is diffeomorphism from $\ID$ to $f(\ID)$ and thus also a minimiser in $\mathsf{F}_p$ should it be continuous on $\overline{\ID}$. We start with the definition of the enlarged space - in fact we will start with the inverse functions.\\

Let $1<p<\infty$. Let $f_0$ be the given boundary data with $\mathsf{E}_p(f_0)<\infty$, and $h_0=f_0^{-1}$. We set $\mathsf{H}_p$ as the space of functions $h:\overline{\ID}\rightarrow\overline{\ID}$ satisfying the following three conditions:\\

\begin{itemize} 
\item $h\in C(\overline{\ID})$, $h|_{\partial\ID}=h_0|_{\partial\ID}$, $h$ has finite distortion, and
\begin{equation}\label{5.1}
\mathsf{E}_p^\ast(h):=\int_\ID\mathbb{K}^p(w,h)J(w,h)dw\leq\int_\ID\mathbb{K}^p(w,h_0)J(w,h_0)dw+1.
\end{equation}
\item Let $q=\frac{2p}{p+1}$. There is an $f\in W^{1,q}(\ID,\ID)$ such that
\begin{equation}\label{5.2}
\|Df\|_{L^q(\ID)}\leq\Big(\pi^p\int_\ID\mathbb{K}^p(z,f_0)dz\Big)^{\frac{1}{p+1}}+1,
\end{equation}
\item There is a measurable set $X\subset\ID$ such that $|\ID-X|=0$, $h\circ f(z)=z$ for every $z\in X$, and $J(w,h)=0$ for almost every $w\in\ID-f(X)$.
\end{itemize}

With the help of (\ref{1.7})-(\ref{1.8}) it is not hard to establish the following lemma.
\begin{lemma}
Every homeomorphic $f\in\mathsf{F}_p$ with
\[
\int_\ID\mathbb{K}^p(z,f)dz\leq\int_\ID\mathbb{K}^p(z,f_0)dz
\]
has its inverse $h=f^{-1}\in\mathsf{H}_p$. In particular,
\[
\inf_{h\in\mathsf{H}_p}\mathsf{E}_p^\ast(h)\leq\inf_{f\in\mathsf{F}_p}\mathsf{E}_p(f).
\]
where $\mathsf{E}_p^\ast$ is defined at (\ref{5.1}).
\end{lemma}

Note that $h_0\in\mathsf{H}_p$, so at least $\mathsf{H}_p\neq\emptyset$. Now let $h_j$ be any sequence in $\mathsf{H}_p$. By (\ref{5.1}), $h_j$ has a uniform $W^{1,2}(\ID)$ norm. So up to a subsequence there is an $h$ such that $h_j\rightharpoonup h$ in $W^{1,2}(\ID)$. Our first task is to show that $\mathsf{H}_p$ is closed under this weak convergence. That is
\begin{theorem}\label{thm 5.1}
Let $h_j\in\mathsf{H}_p$, $h_j\rightharpoonup h$ in $W^{1,2}(\ID)$. Then $h\in\mathsf{H}_p$, and
\[
\mathsf{E}_p^\ast(h)\leq\liminf_{j\rightarrow\infty}\mathsf{E}_p^\ast(h_j).
\]
\end{theorem}

We first remark that the local convergence can be extended to $\overline{\ID}$. In fact, each $h_j$ can be continuously extended to some $\ID_R$ ($R>1$) with the same function $\frac{1}{\overline{h_0(\frac{1}{\bar{w}})}}$, $w\in\ID_R\setminus\ID$. Then the local uniform convergence applies on $\ID_R$. Precisely,

\begin{lemma}
Let $h_j\in\mathsf{H}_p$ be a sequence such that $h_j\rightharpoonup h$ weakly in $W^{1,2}(\ID)$. Then there is a subsequence $h_j\rightrightarrows h$ uniformly in $\overline{\ID}$ and $J(w,h_j)\rightharpoonup J(w,h)$ weakly in $L^1(\ID)$.
\end{lemma}

\begin{lemma}
$h$ is a $W^{1,2}(\ID)\cap C(\overline{\ID})$ finite distortion function, $h|_{\partial\ID}=h_0|_{\partial\ID}$, and it satisfies (5.1).
\end{lemma}
\noindent{\bf Proof.} As a weak limit is clear that $h\in W^{1,2}(\ID)$. Since $h_j\rightrightarrows h$ uniformly on $\overline{\ID}$, we also have that $h\in C(\overline{\ID})$ with the same boundary values $h_0|_{\partial\ID}$. For the finite distortion, we have the following inequality [9, Lemma 8.8.2]:
\[
x^ny^{-l}-x_0^ny_0^{-l}\geq nx_0^{n-1}y_0^{-l}(x-x_0)-lx_0^ny_0^{-l-1}(y-y_0),
\]
for $n\geq l+1\geq1$. Put $\|Dh_j(w)\|$, $J(w,h^j)$, $\|Dh(w)\|$ and $J(w,h)$ into it we get
\begin{align*}
&\frac{\|Dh_j(w)\|^{2p}}{J^{p-1}(w,h_j)}-\frac{\|Dh(w)\|^{2p}}{J^{p-1}(w,h)}\notag\\
\geq&2p\frac{\|Dh(w)\|^{2p-1}}{J^p(w,h)}(\|Dh_j\|-\|Dh\|)-(p-1)\frac{\|Dh(w)\|^{2p}}{J^p(w,h)}(J(w,h_j)-J(w,h)).
\end{align*}
Upon integration the right-hand side here converges to $0$, c.f. \cite[Theorem 12.2]{AIMO}. Then it follows that
\begin{equation}\label{5.3}
\mathsf{E}_p^\ast(h)\leq\liminf_{j\rightarrow\infty}\mathsf{E}_p^\ast(h_j).
\end{equation}
So $h$ satisfies (\ref{5.1}) as each $h_j$ does, and in particular, $h$ has finite distortion.\hfill $\Box$

\bigskip

\begin{lemma}
Every $W^{1,2}(\ID)$ finite distortion function satisfies Lusin's condition $\mathcal{N}$.
\end{lemma}
This is proved by Gol'dshtein and Vodop'yanov in \cite{GV}.\\

Now we write $f_j$, $X_j$ as in the definition of each $h_j$.  Directly from the definition,  $f_j$ are bounded in $L^q(\ID)$, so there is a subsequence such that $f_j\rightharpoonup f$ in $W^{1,q}(\ID)$. By the lower semi-continuity of weak convergence we have
\[
\|Df\|_{L^q(\ID)}\leq\liminf_{j\rightarrow\infty}\|Df_j\|_{L^q(\ID)}\leq\Big(\pi^p\int_\ID\mathbb{K}^p(z,f_0)dz\Big)^{\frac{1}{p+1}}+1.
\]
So $f$ satisfies (\ref{5.2}). On the other side, by the Rellich-Kondrachov Theorem, we have $f_j\rightarrow f$ strongly in $L^s(\ID)$, for all $1\leq s<q^*$, where $q^*=\frac{2q}{2-q}$ is the Sobolev conjugate. In particular, again up to a subsequence we have that $f_j\rightarrow f$ pointwise almost everywhere in $\ID$.\\

\begin{lemma}
Let $g\in\mathsf{H}_p$ and $f_g$, $X_g$ be the corresponding function and set as in the definition of the space $\mathsf{H}_p$. Then, for every measurable function $\eta$ defined on $\ID$,
\[
\int_\ID\eta(w)J(w,g)dw=\int_\ID\eta(f_g(z))dz.
\]
\end{lemma}
\noindent{\bf Proof.} From the assumptions  $|\ID-X_g|=0$ and $J(w,g)=0$ almost everywhere  in $\ID-f_g(X_g)$, the equality reads as
\[
\int_{f_g(X_g)}\eta(w)J(w,g)dw=\int_{X_g}\eta(f_g(z))dz.
\]
However, this is simply the area formula together with Lusin's condition $\mathcal{N}$ for $g$. See \cite[Theorem 2]{H}.\hfill $\Box$

\bigskip

Now define
\[
X_h:=\{z\in\ID:f_j(z)\rightarrow f(z)\}\cap \bigcap_{j=1}^\infty X_j.
\]
Note $X_h$ still has full measure in $\ID$.

\begin{lemma}
For every $z\in X_h$, $h\circ f(z)=z$.
\end{lemma}
\noindent{\bf Proof.}
\begin{align*}
|z-h(f(z))|&=|h_j(f_j(z))-h(f(z))|\\
&\leq|h_j(f_j(z))-h(f_j(z))|+|h(f_j(z))-h(f(z))| \rightarrow0,
\end{align*}
as $h_j\rightrightarrows h$ in $\overline{\ID}$.\hfill $\Box$

\medskip

\begin{lemma}
\[
J(w,h)=0,\quad a.e.\quad w\in\ID-f(X_h).
\]
\end{lemma}
\noindent{\bf Proof.} Let $\eta\in C(\ID)\cap L^\infty(\ID)$. Then,
\begin{eqnarray*}
\lefteqn{\int_{\ID}\eta(w)J(w,h)dw}\\ &=&\lim_{j\rightarrow\infty}\int_{\ID}\eta(w)J(w,h_j)dw  \lim_{j\rightarrow\infty}\int_\ID\eta(f_j(z))dz = \int_\ID\eta(f(z))dz,
\end{eqnarray*}
since $J(w,h_j)\rightharpoonup J(w,h)$ in $L^1(\ID)$ and $f^j\rightarrow f$ pointwise almost everywhere  in $\ID$. 

We now let $\eta^k\rightarrow\chi_{\ID-f(X)}$ be the standard mollification. We only need the pointwise convergence and the fact $\|\eta^k\|_\infty\leq\|\chi_{\ID-f(X)}\|_\infty\leq1$, which is a property of convolutions. Now by dominated convergence,
\begin{eqnarray*}
\lefteqn{\int_{\ID}\chi_{\ID-f(X_h)}(w)J(w,h)dw}\\&=&\lim_{k\rightarrow\infty}\int_{\ID}\eta^k(w)J(w,h)dw = \lim_{k\rightarrow\infty}\int_\ID\eta^k(f(z))dz = \int_\ID\chi_{\ID-f(X_h)}(f(z))dz=0.
\end{eqnarray*}
Note here $\eta^k(f(z))\rightarrow\chi_{\ID-f(X)}(f(z))$ pointwise almost everywhere  in $\ID$ because $h$ satisfies Lusin's condition $\mathcal{N}$.\hfill $\Box$

\medskip

We have now verified that $h$, $f$ and $X_h$ satisfy all the conditions of $\mathsf{H}_p$, so Theorem \ref{thm 5.1} is proved.\\

We now let $h_j\in\mathsf{H}_p$ be a minimising sequence. Then, there is a limit function $h\in\mathsf{H}_p$. By (\ref{5.3}), $h$ is a minimiser. To see $h$ is variational, we need the following lemma.

\begin{lemma} [Chain rule]
Let $f\in W^{1,1}(\Omega,\Omega^\prime)$, $h\in W^{1,1}(\Omega^\prime,\IC)\cap C(\Omega^\prime)$. Assume that $h\circ f\in W^{1,1}(\Omega)$, and $f$ has Lusin $\mathcal{N}^{-1}$. Then, for almost every $z\in\Omega$,
\begin{equation}\label{5.4}
(h\circ f)_z(z)=h_w(f(z))f_z(z)+h_{\bar{w}}(f(z))\overline{f_{\bar{z}}}(z),
\end{equation}
\begin{equation}\label{5.5}
(h\circ f)_{\bar{z}}(z)=h_w(f(z))f_{\bar{z}}(z)+h_{\bar{w}}(f(z))\overline{f_z}(z).
\end{equation}
\end{lemma}
\noindent{\bf Proof.} We prove (\ref{5.4}). Let $h^\varepsilon\to h$ be the standard mollification. Then
\[
(h^\varepsilon\circ f)_z(z)=h^\varepsilon_w(f(z))f_z(z)+h^\varepsilon_{\bar{w}}(f(z))\overline{f_{\bar{z}}}(z).
\]
It is clear that for almost every point $z$, the right hand side converges to
\[
h_w(f(z))f_z(z)+h_{\bar{w}}(f(z))\overline{f_{\bar{z}}}(z),
\]
and so we need to show
\[
\lim_{\varepsilon\to0}(h^\varepsilon\circ f)_z(z)=(h\circ f)_z(z),\quad a.e.\quad z\in\ID.
\]
We let $\eta$ be the standard mollifier. Then
\begin{eqnarray}
\lefteqn{\lim_{\varepsilon\to0}(h^\varepsilon\circ f)_z(z)}\notag\\&=&\lim_{\varepsilon\to0}\lim_{k\to0}\int_\ID(h^\varepsilon\circ f)_z(\zeta)\eta^k(z-\zeta)d\zeta
=\lim_{\varepsilon\to0}\lim_{k\to0}-\int_\ID h^\varepsilon\circ f(\zeta)\eta^k_z(z-\zeta)d\zeta\notag\\
&=&\lim_{k\to0}\lim_{\varepsilon\to0}-\int_\ID h^\varepsilon\circ f(\zeta)\eta^k_z(z-\zeta)d\zeta \label{5.6} = \lim_{k\to0}-\int_\ID h\circ f(\zeta)\eta^k_z(z-\zeta)d\zeta \\
&=&\lim_{k\to0}\int_\ID (h\circ f)_z(\zeta)\eta^k(z-\zeta)d\zeta= (h\circ f)_z(z).\notag
\end{eqnarray}
We explain the interchange of these two limits in (\ref{5.6}). In fact, this holds if one of them is uniform. However, we know $h^\varepsilon\to h$ uniformly in a neighbourhood of $z$, since $h$ is continuous. Thus, for any $\delta$, we may choose $\varepsilon_0$ so small that for any $\varepsilon<\varepsilon_0$, $\|h^\varepsilon-h\|_\infty<\delta$. Then, for any $k$ 
\begin{eqnarray*}
\lefteqn{\Big|\int_\ID h^\varepsilon\circ f(\zeta)\eta^k_z(z-\zeta)d\zeta-\int_\ID h\circ f(\zeta)\eta^k_z(z-\zeta)d\zeta\Big|}\\
&\leq&\|h^\varepsilon-h\|_\infty\int_\ID|\eta^k_z(z-\zeta)|d\zeta \leq \delta \frac{1}{k^2}\int_{B(z,k)}|\eta_z(\frac{z-\zeta}{k})|d\zeta\\
&\leq&\delta\pi\|\nabla\eta\|_\infty
\end{eqnarray*}
This proves the uniform convergence. \hfill$\Box$

\bigskip

Now our $h$ and $f$ satisfy $h\circ f(z)=z$ almost everywhere  $z\in\ID$, so by (\ref{5.4}) and (\ref{5.5}),
\[
f_z(z)=\frac{\overline{h_w(f(z))}}{J(f(z),h)},\quad f_\zbar(z)=-\frac{h_\wbar(f(z))}{J(f(z),h)},\quad\mbox{a.e. }z\in\ID.
\]
In particular,
\[
\IK(f(z),h)=\IK(z,f),\quad J(f(z),h)J(z,f)=1,\quad\mbox{a.e. }z\in\ID.
\]
Then, as before,
\begin{align*}
\mathsf{E}_p^\ast(h)=\int_\ID\mathbb{K}^p(w,h)J(w,h)dw& 
 =\int_\ID\mathbb{K}^p(z,f)dz=\mathsf{E}_p(f).
\end{align*}
Then, when $h_j$ is a minimising sequence for $\mathsf{E}_p^\ast$, $f_j$ is also a minimising sequence for $\mathsf{E}_p$. This implies that
\[
\int_\ID\mathbb{K}^p(w,h)J(w,h)dw\leq\int_\ID\mathbb{K}^p(w,h_0)J(w,h_0)dw<\int_\ID\mathbb{K}^p(w,h_0)J(w,h_0)dw+1,
\]
and
\begin{align*}
\|Df\|_{L^q(\ID)}&\leq\Big(\pi^p\int_\ID\mathbb{K}^p(z,f)dz\Big)^{\frac{1}{p+1}}\\
&\leq\Big(\pi^p\int_\ID\mathbb{K}^p(z,f_0)dz\Big)^{\frac{1}{p+1}}<\Big(\pi^p\int_\ID\mathbb{K}^p(z,f_0)dz\Big)^{\frac{1}{p+1}}+1.
\end{align*}
This implies $h$ is inner-variational in the space $\mathsf{H}_p$, and then $h$ and $f$ satisfy the variational equations (\ref{1.5}) and (\ref{4.1}). In particular, if $f$ satisfies either of the conditions  (1)  to  (6)  in Theorem \ref{4.1}, then  $f$ is a diffeomorphism. Furthermore, since $\{h:h^{-1}\in\mathsf{F}_p\}\subset\mathsf{H}_p$, $f$ is also a minimiser in $\mathsf{F}_p$. We have proved the following.
\begin{theorem}\label{thm 5.2}
The space $\mathsf{H}_p$ admits a minimiser $h$. Let $f$ be its `inverse' as in the definition of $\mathsf{H}_p$. Then $h$ and $f$ satisfies the inner-variational equations
\[
2p\int_\ID \IK_f^p\frac{\overline{\mu_f}}{1+|\mu_f|^2}\varphi_\zbar dz=\int_\ID \IK_f^p\varphi_zdz,\quad\forall\varphi\in C_0^\infty(\ID ),
\]
and
\[
\int_\ID \IK^{p-1}(w,h)h_w\overline{h_\wbar }\varphi_\wbar dw=0,\quad\forall\varphi\in C_0^\infty(\ID ).
\]
Furthermore, if any of the conditions(1) to (6) in Theorem \ref{4.1} is satisfied, then $f$ is a diffeomorphic minimiser in $\mathsf{F}_p$.
\end{theorem}

\section{A topological condition.} 

In this section we show that if the minimiser   of Theorem \ref{thm 5.2} lies in the Sobolev space $W^{1,2}(\ID)$,  then it is in the space $\mathsf{F}_p$.  That is $f$ is a homeomorphism and $f|\IS=f_0|\IS$.  We start with a theorem giving the topological result under an assumption about the existence of a {\em principal solution}.  These solutions are discussed more fully in \cite[\S 20.2]{AIM}.  We only need here that they are entire $W^{1,2}_{loc}(\IC)$ homeomorphic solutions to a Beltrami equation normalised so as to be conformal near $\infty$.

\begin{theorem}\label{thm 6.1} Let $f:\ID\to\ID$ be a surjective mapping of finite distortion and topological degree $1$. Suppose $f\in W^{1,2}(\ID)$ and $\IK(z,f)\in L^1(\ID)$ and also that there is a $W^{1,2}_{loc}(\IC)$ homeomorphic solution $F:\IC\to\IC$ to the Beltrami equation
\begin{equation}
F_\zbar=\mu(z)\; F_z,\hskip15pt a.e. \;\; z\in \IC.
\end{equation}
with $\mu(z)=\mu_f(z)$ for $z\in \ID$.
Then $f:\overline{\ID}\to\overline{\ID}$ is a homeomorphism.
\end{theorem}
\noindent{\bf Proof.}  Let $\Omega=F(\ID)$.  Then as $F$ is a homeomorphism of $\IC$ we see $\Omega$ is a Jordan domain.  The mapping $H=F^{-1}:\Omega \to \ID$ has 
\[ \int_\Omega \|\nabla H\|^2 \; dw = \int_\ID \|\nabla H(F)\|^2 J(z,F) \; dz = \int_\ID \frac{\|\nabla F\|^2 }{J(z,F)}\; dz  = \int_\ID \IK(z,F) \; dz <\infty \]
since under these hypothesis,  as noted above,  the change of variables formula holds.  We may therefore calculate that
\begin{eqnarray*}
\lefteqn{\int_\Omega \| \nabla ( f \circ H) \| \; dw } \\ &=& \int_\ID \| \nabla f (z) \nabla H(F) \| \; J(z,F) \; dz   \leq   \int_\ID \| \nabla f (z)\|  {\| (\nabla F)^{-1} \|}\,{{J(z,F)}} \; dz \\
& = &     \int_\ID \| \nabla f (z)\|  \, \| \nabla F (z) \|   \; dz  \leq   \Big( \int_\ID \| \nabla f \|^2 \; dz \Big)^{\frac{1}{2}} \; \Big(\int_\ID \|\nabla F\|^2  \; dz \Big)^{\frac{1}{2}} 
\end{eqnarray*}
Thus $f\circ H \in W^{1,1}(\Omega)$.  The degree of this bounded mapping is $1$.  We next calculate that
\begin{equation}
|\mu_{f\circ H}(F(z))| = \left| \frac{\mu_f(z)-\mu_F(z)}{1-\bar\mu_f(z) \mu_F(z)}\right| = 0, \hskip10pt a.e.\;\; z\in \Omega
\end{equation}
In view of Weyl's lemma,  we have now shown $\phi=f\circ H:\Omega\to\ID$ is onto, holomorphic and degree $1$ and therefore is a conformal mapping.  Since $\Omega$ is a Jordan domain $\phi$ extends homeomorphically to the boundary of $\Omega$ by Carath\'eodory's theorem.  We have now that $\phi(F)=f:\ID\to\ID$ and the left-hand side extends homeomorphically to the boundary.  Therefore the right-hand side does as well and this proves the theorem.\hfill $\Box$

\medskip

\begin{theorem}\label{thm 6.2} Let $ h \in \mathsf{H}_p$ with $\IK(w,h)\in L^1(\ID)$ and quasiconformal boundary data $h_0:\overline{\ID}\to\overline{\ID}$.  Then $h$ is a homeomorphism.
\end{theorem}
\noindent{\bf Proof.} In order to apply Theorem \ref{thm 6.1} we seek a principal solution to the equation 
\begin{equation}
H_\wbar=\mu\; H_w,\hskip15pt a.e. \;\; w\in \IC.
\end{equation}
for any $\mu$ equal to $\mu_h$ on $\ID$.  Now,  that $h\in \mathsf{H}_p$ gives us $h\in W^{1,2}(\ID)$ is continuous on $\overline{\ID}$ and $h|\IS=h_0$.  Define a mapping
\begin{equation}
h^\ast = \left\{\begin{array}{ll} h(w), & w\in \ID,  \\ 1/\overline{ h_0(\frac{1}{\wbar})}, & w \in \IC\setminus \ID. \end{array}\right.
\end{equation}
The mapping $h^\ast$ is continuous and quasiconformal in $\IC\setminus \ID$. Suppose $h_0(w_0)=0$ and set $r=(1+|w_0|)/2$. Let $R=2\max_{z\in \IS(r)} |h^\ast(z)| <\infty$. Then $h^\ast|\IS(r)$ and $identity|\IC\setminus\ID(R)$ are a pair of quasiconformal embeddings of $\IS(r)$ and $\IS(R)$ with disjoint images.  The quasiconformal version of the Sch\"onflies Theorem, \cite{GVam} or \cite[\S 7]{GMP} tells us that there is a quasiconformal mapping $g:A(r,R)=\{r<|z|<R\} \to \IC$ with $g|\IS(r)=h^\ast$ and $g|\IS(R)=identity$.  Now define a mapping $g^\ast:\IC\to\IC$ by
\begin{equation}
g^\ast = \left\{\begin{array}{ll} h^\ast(w), & w\in \overline{\ID(r)},  \\  g(w), & w \in \overline{A(r,R)} \\ w, & |w|\geq R \end{array}\right.
\end{equation}
The mapping $g^\ast$ is quasiconformal on $\IC\setminus \ID$,  conformal on $\IC\setminus\ID(R+1)$.  It follows that $g^\ast\in W^{1,2}(\ID(R+1))$ and that $\mu_{g^\ast}=0$ outside of $\ID(R)$.  Now  \cite[Theorem 20.2.1]{AIM} provides a principle solution (we have to make the minor adjustment of replacing $\ID$ by $\ID(0,R+1)$on which $\mu$ is compactly supported).  The result follows. \hfill $\Box$

\medskip

\noindent{\bf Remark 1.} Here the boundary values $h_0:\ID\to\ID$ do not need to be quasiconformal,  though some restriction is required.   It is easy to see that the proof given works as soon as the boundary values $h_0$ admit an extension to an annulus $A(r,1)$ with both $h_0$ and $(h_0)^{-1}$ having finite Dirichlet energy and $h_0$ locally quasiconformal on a neighbourhood of some $\IS(s)$, $r<s<1$.  We then simply use the Sch\"onflies Theorem on $\IS(1/s)$ via the reflection of $h_0$. The local quasiconformality would be implied,  for instance, by the extension being a diffeomorphism on $A(r,1)$.  We identify a necessary condition below.

\begin{theorem}. Let $f_0:\IS\to\IS$ be a homeomorphism.  Then $f_0$ admits an extension $f:\overline{\ID}\to\overline{\ID}$ with $f|\IS=f_0$,  $f\in W^{1,2}(\ID)$ and $\IK(z,f)\in L^1(\ID)$ only if 
\begin{equation}
\iint_{\IS\times\IS}\big( Q + \big|\,\log Q\,\big| \big) \; |d\zeta||d\xi| < \infty,  \hskip10pt Q(\zeta,\xi) = \left|\frac{f_0(\zeta)-f_0(\xi)}{\zeta-\xi}\right|^2.
\end{equation}
\end{theorem}
\noindent{\bf Proof.} The condition $\iint_{\IS\times\IS}  Q |d\zeta||d\xi|<\infty$ is Douglas' necessary and sufficient condition for $f_0$ to admit an extension of finite Dirichlet energy, \cite{Douglas}.  In \cite{AIMO} the condition 
\begin{equation}\label{6.7}   \iint_{\IS\times\IS} - \log |f_0(\zeta)-f_0(\xi) |  \; d\zeta d\bar \xi   < \infty \end{equation}
 is shown to be necessary and sufficient for $f_0$ to admit an extension with $\IK(z,f)\in L^1(\ID)$.  Since
\[  \iint_{\IS\times\IS}  \log |\zeta-\xi | d\zeta d\bar \xi = - \pi \]
and that this integral is uniformly convergent,  we see that the integral at (\ref{6.7}) is finite if and only if 
\[ \iint_{\IS\times\IS}   \big| \log Q \big| \; |d\zeta|| d\bar \xi |<\infty \]
and this completes the proof. \hfill $\Box$

\medskip

\noindent{\bf Remark 2.}  The hypothesis $\IK(w,h)\in L^1(\ID)$ is equivalent to the condition $f\in W^{1,2}(\ID)$ for the pseudo-inverse $f=h^{-1}$ found in the definition of the space $\mathsf{H}_p$.  We call the pseudo-inverse $f$ of $h$  a minimiser of $\mathsf{E}_p$ when $h$ is a minimiser as described by Theorem \ref{thm 5.2}.

\medskip

In any case, with Theorems \ref{thm 4.1} and \ref{thm 6.2} we have the following result.

\begin{theorem}\label{6.3} Let $p>1$ and $f_0:\IS\to\IS$ be quasisymmetric boundary data.  Let $f=h^{-1}$ where $h \in  \mathsf{H}_p$ is the minimiser provided by Theorem \ref{thm 5.2}. Then $f\in W^{1,\frac{2p}{p+1}}(\ID)$,  and
\begin{enumerate}
\item  $f\in W^{1,2}(\ID)$ implies  $f:\overline{\ID}\to\overline{\ID}$ is a homeomorphism,  while
\item  $f\in W^{1,2}(\ID)\cap W^{1,s}_{loc}(\ID)$ implies $f:\ID\to\ID$ is also a diffeomorphism.
\end{enumerate}
\end{theorem}

\section{Limiting regimes.}

In this section we discuss what happens as $p\to\infty$ or $p\to 1$ for fixed boundary data $f_0:\IS\to\IS$.  Theorem \ref{6.3} suggests improved regularity as $p\to\infty$ for then $\frac{2p}{p+1}\to 2$,  while as $p\to 1$ we only have the weaker bounds as $\frac{2p}{p+1}\to 1$.  It is a little surprising then that the minimisers of $\mathsf{E}_\infty$ are always (quasiconformal) homeomorphisms,  but not necessarily diffeomorphisms,  while for $\mathsf{E}_1$ minimisers are always diffeomorphisms,  but almost never quasiconformal.

\subsection{$p\to \infty$}

We first make a definition.  We say a sequence of holomorphic functions $\{\psi_k\}$ with $\psi_k\in L^1(\ID)$  is degenerate if  
\begin{equation}
\frac{\psi_k}{\|\psi_k\|_{L^1(\ID)}} \to 0, \hskip10pt \mbox{locally uniformly in $\ID$}.
\end{equation}
Otherwise the sequence is nondegenerate.

\begin{theorem} \label{thm 7.1} Let $h_0:\overline{\ID }\to\overline{\ID }$ be quasisymmetric For each $p$ let $h_p\in\mathsf{F}_p$ be minimiser for the boundary values $h_0$, and $\phi_p$ the associated Ahlfors-Hopf holomorphic quadratic differential.  Then the following hold.
\begin{enumerate}
\item There is a quasiconformal $h:\ID\to\ID$ and a subsequence $\{p_k\}_{k=1}^{\infty}$ so that
$h_{p_k}\to h$ uniformly in $\ID$ and weakly in $W^{1,2}(\ID)$.
\item The mapping $h$ is extremal for its boundary values and
\[ \liminf [\mathsf{E}_p^{\ast}(h_p)]^{1/p} = \| \IK(w,h)\|_{L^\infty(\ID)} \]
\item If the sequence $\{ \phi_{p_k} \}$ is nondegenerate, then the mapping $h$ is a Teichm\"uller mapping,  $\mu_h(w) = k \frac{\overline{\psi}}{|\psi|}$, $k\in [0,1)$ and $\phi\in L^1(\ID)$ holomorphic. Then $h$ is uniquely extremal and $h_p\to h$ uniformly.
\end{enumerate}
\end{theorem}

The hypothesis of nondegeneracy in {\em 3.} above is necessary.  There are uniquely extremal quasiconformal self mappings of $\ID$ with nonconstant distortion $\IK(z,h)$, see Mateljevi\'c's survey \cite[pp 86--88]{Mat} and the references therein.  For such a mapping {\em 1.} and {\em 2.} show that $h_p\to h$ uniformly in $\ID$,  while {\em 3.} shows the sequence of Ahlfors-Hopf differentials must be degenerate.

\medskip

This theorem has a couple of corollaries which follow from its proof.  In view of the distributional equation one might conjecture that for an extremal $\IK(z,f)\in L^p(\ID)$ should not have $\mu=0$ (consider test functions supported near this zero).  However this not not the case.  First a smoothness condition.

 \begin{corollary}\label{cor1} With the notation of Theorem \ref{thm 7.1},  suppose that $p_k\to\infty$, $\{\phi_{p_k} \}$ is nondegenerate and  $\phi_{p_k} \neq 0$ on $\ID$.  Then the unique extremal quasiconformal mapping $h$ is a diffeomorphism of $\ID$.
 \end{corollary}
 
 \begin{corollary}\label{cor2} There are quasisymmetric mappings $f_0:\IS\to\IS$ for which for all sufficiently large $p$,  the extremal $h_p\in\mathsf{F}_p$ has a point $w_0\in \ID$ where $\mu_f(w_0)=0$ and $\phi_p(w_0)=0$.  \end{corollary}

Following Ahlfors,  we give a proof by considering $p\to\infty$ in the $L^p$ problems, and point out where the differences lie. 
We recall that
\begin{equation}\label{7.2}
\phi_p=\IK_{h_p}^{p-1}(h_p)_w\overline{(h_p)_\wbar }=\IK_{h_p}^pJ_{h_p}\frac{\overline{\mu_{h_p}}}{1+|\mu_{h_p}|^2},
\end{equation}
and also that the sequence is uniformly bounded in $W^{1,2}(U)$ where $U$ is open and $\overline{\ID}\subset U$.  Thus let $h$ be the weak limit of $h_p$ in $U$ and the uniform convergence on $\overline{\ID}$ is assured,  $h\in W^{1,2}(\ID)$ is continuous and $h|\IS=h_0$.

First of all, if $\frac{1}{\pi}\mathsf{E}_\infty^\ast(h)=1$, then $h$ is a conformal mapping as $\IK(z,f)\equiv 1$ and $k=0$.  Henceforth, we assume that $\frac{1}{\pi}\mathsf{E}_\infty^\ast(h)>1$. Set
\begin{equation}\label{7.2}
k:=\sqrt{\frac{\frac{1}{\pi}\mathsf{E}_\infty^\ast(h)-1}{\frac{1}{\pi}\mathsf{E}_\infty^\ast(h)+1}}>0.
\end{equation}
\begin{lemma}\label{lem 7.1}
\begin{equation}
\liminf_{p\to\infty}C_p^\frac{1}{p}=\mathsf{E}_\infty^\ast(h), \hskip10pt C_p:=\int_\ID |\phi_p|.
\end{equation}
\end{lemma}
\noindent{\bf Proof.} We first observe
\begin{align*}
C_p=\int_\ID\IK_{h_p}^pJ_{h_p}\frac{|\mu_{h_p}|}{1+|\mu_{h_p}|^2}&\leq\frac{1}{2}\int_\ID\IK_{h_p}^pJ_{h_p} =\frac{1}{2}\mathsf{E}_p^{\ast}(h_p)=\frac{1}{2}\mathsf{E}_p(f_p).
\end{align*}
Thus
\[
\lim_{p\to\infty}\Big(\int_\ID |\Psi_p(w)|\Big)^\frac{1}{p}\leq\lim_{p\to\infty} \frac{1}{2^\frac{1}{p}}\mathsf{E}_p^\frac{1}{p}(f_p)=\mathsf{E}_\infty(f_\infty).
\]
This proves one direction. For the other direction, observe
\[
\IK_{h_p}^p\leq\IK_{h_p}^p\frac{|\mu_{h_p}|}{1+|\mu_{h_p}|^2}\frac{1+\delta_p^2}{\delta_p}+\Big(\frac{1+\delta_p^2}{1-\delta_p^2}\Big)^p,\quad\forall\delta_p\in(0,1).
\]
Multiply by $J_{h_p}$ and integrate both sides over $\ID $ to obtain
\[
\mathsf{E}_p(f_p)=\mathsf{E}_p^\ast(h_p)\leq\frac{1+\delta_p^2}{\delta_p}\int_\ID |\Psi_p(w)|dw+\pi\Big(\frac{1+\delta_p^2}{1-\delta_p^2}\Big)^p.
\]
Now,  for each $p$ there is a $\delta_p\in(0,1)$ such that
\[
\Big(\frac{1+\delta_p^2}{1-\delta_p^2}\Big)^p=\frac{1}{2}\mathsf{E}_p(f_p).
\]
Then
\[
\frac{1+\delta_p^2}{1-\delta_p^2}=\Big(\frac{1}{2}\mathsf{E}_p(f_p)\Big)^\frac{1}{p}\to\mathsf{E}_\infty(f_\infty).
\]
That is,
\[
\lim_{p\to\infty}\delta_p=\sqrt{\frac{\mathsf{E}_\infty(f_\infty)-1}{\mathsf{E}_\infty(f_\infty)+1}}=k.
\]
So for every $p$,
\[
\int_\ID |\Psi_p(w)|\geq\frac{\delta_p}{1+\delta_p^2}\cdot\frac{1}{2}\mathsf{E}_p(f_p)\geq\frac{\delta_p}{4}\mathsf{E}_p(f_p),
\]
and then
\[
\lim_{p\to\infty}\Big(\int_\ID |\Psi_p(w)|\Big)^\frac{1}{p}\geq\lim_{p\to\infty}\big(\frac{\delta_p}{4}\big)^\frac{1}{p}\mathsf{E}_p^\frac{1}{p}(f_p)=\mathsf{E}_\infty(f_\infty).
\]

\hfill $\Box$

\medskip

Next, with $k$ defined at (\ref{7.2}) we want to consider the possibility that
\begin{equation}\label{7.5}
\lim_{p\to\infty}\int_\ID \Big||(h_p)_\wbar |-k|(h_p)_w|\Big|\to0.
\end{equation}
To analyse this limit choose any $\varepsilon>0$, so small that both $k(1+\varepsilon)$ and $k(1-\varepsilon)$ are in $(0,1)$, and define the two sets
\[
E_p:=\{w\in\ID :|\mu_{h_p}(w)|>k(1+\varepsilon)\},
\]
\[
F_p:=\{w\in\ID :|\mu_{h_p}(w)|<k(1-\varepsilon)\}.
\]
First,
\begin{eqnarray*}
C_p&=&\int_\ID\IK_{h_p}^{p-1}|(h_p)_w(h_p)_\wbar| \geq \int_{E_p}\IK_{h_p}^{p-1}|(h_p)_w(h_p)_\wbar|\\
&\geq&\Big(\frac{1+k^2(1+\varepsilon)^2}{1-k^2(1+\varepsilon)^2}\Big)^{p-1}k(1+\varepsilon)\int_{E_p}|(h_p)_w|^2.
\end{eqnarray*}
Rearrange this to obtain
\begin{equation}
\int_{E_p}|(h_p)_w|^2\leq\frac{1+k^2(1-\varepsilon)^2}{k(1+\varepsilon)[1+k^2(1+\varepsilon)^2]}\Big(C_p^\frac{1}{p}\frac{1-k^2(1+\varepsilon)^2}{1+k^2(1+\varepsilon)^2}\Big)^p.
\end{equation}
By Lemma \ref{lem 7.1}, if $|E_p|\neq 0$ we know that for $p$ sufficiently large,
\begin{equation}\label{7.6}
\liminf C_p^\frac{1}{p}  >  \frac{1+k^2(1+\varepsilon)^2}{1-k^2(1+\varepsilon)^2},
\end{equation}
and this implies $
\lim_{p\to\infty}\int_{E_p}|(h_p)_w|^2=0
$ 
and hence
\[
\lim_{p\to\infty}\int_{E_p}(|(h_p)_\wbar |-k|(h_p)_w|)^2\leq(1+k)^2\lim_{p\to\infty}\int_{E_p}|(h_p)_w|^2=0.
\]
Now applying the Cauchy-Schwarz inequality gives
\[
\lim_{p\to\infty}\int_{E_p}\Big||(h_p)_\wbar |-k|(h_p)_w|\Big|=0.
\]
At this point we have established both {\em 1.} and {\em 2.} of Theorem \ref{thm 7.1}.  Thus we now turn to the set $F_p$. Set

\begin{equation}
\Xi_p=\frac{\phi_p}{\|\phi_p\|_{L^1(\ID )}}=\frac{\phi_p}{C_p}.
\end{equation}
Then the sequence  $\{\Xi_p,p\geq1\}$ is bounded in $L^{1}(\ID)$ and so form a normal family of analytic functions there, \cite{HM}.  
We find a subsequence converging locally uniformly to
\begin{equation}\label{7.8}
\psi=\lim_{p_k\to\infty}\Xi_p.
\end{equation} 
We compute
\begin{eqnarray*}
\int_{F_p}|\Xi_p|&=&\frac{1}{C_p}\int_{F_p}\IK_{h_p}^{p-1}|(h_p)_w(h_p)_\wbar|\\
&\leq & \frac{1}{C_p}\Big(\frac{1+k^2(1-\varepsilon)^2}{1-k^2(1-\varepsilon)^2}\Big)^{p-1}k(1-\varepsilon)\int_{F_p}|(h_p)_w|^2
\end{eqnarray*}
and again Lemma \ref{lem 7.1}  and the fact that $\|Dh_p\|_{L^2(\ID )}$ are uniformly bounded gives
\begin{equation}\label{7.10}
\lim_{p\to\infty}\int_{F_p}|\Xi_p|=0.
\end{equation}

Now we use the assumption of nondegeneracy which implies there is a subsequence converging locally unformly in $\ID$ to
\begin{equation}
\psi=\lim_{p_k\to\infty}\Xi_{p_k}.
\end{equation}
and $\psi$ is not identically $0$.

Now since $\Xi_p\not\to 0$, (\ref{7.10}) implies $|F_p|\to 0$. Thus
\begin{equation}
\lim_{p\to\infty}\int_{F_p}\Big||(h_p)_\wbar |-k|(h_p)_w|\Big|=0.
\end{equation}
Next, in $V_p=\ID\setminus \{E_p\cup F_p\}$,
$
k(1-\varepsilon)\leq|\mu_{h_p}|\leq k(1+\varepsilon)
$
and so 
\begin{eqnarray*}
\lefteqn{\lim_{p\to\infty}\int_{V_p}\Big||(h_p)_\wbar |-k|(h_p)_w|\Big| }\\
&=&\lim_{p\to\infty}\int_{V_p}|(h_p)_w|\Big||\mu_{h_p}|-k\Big| \leq k\varepsilon\lim_{p\to\infty}\int_{V_p}|(h_p)_w| 
 \leq \pi^\frac{1}{2}k\varepsilon\|Dh_p\|_{L^2(\ID )}^\frac{1}{2}.
\end{eqnarray*}
Again we have that $\|Dh_p\|_{L^2(\ID )}$ are uniformly bounded. So (\ref{7.5}) follows as $\varepsilon$ can be arbitrarily small.
\hfill $\Box$

\bigskip

\begin{lemma}
Let $\psi$ be the holomorphic limit of $\Xi_p$ as defined in (\ref{7.8}). Then,
\begin{equation}
\mu_h=k\frac{\overline{\Psi}}{|\Psi|}.
\end{equation}
\end{lemma}
\noindent{\bf Proof.}
From (\ref{7.5})
\begin{eqnarray}
\lefteqn{\int_\ID \Big|(h_p)_\wbar \frac{\Xi_p}{|\Xi_p|}-k(h_p)_w\Big|=\int_\ID \Big|(h_p)_\wbar \frac{\psi_p}{|\psi_p|}-k(h_p)_w\Big|}\notag\\&&=\int_\ID \Big|\frac{(h_p)_w|(h_p)_\wbar |}{|(h_p)_w|}-k(h_p)_w\Big|=\int_\ID \Big||(h_p)_\wbar |-k|(h_p)_w|\Big|\to0.
\end{eqnarray}
We set
\[
\ID _\varepsilon:=\{w\in\ID :|\psi(w)|>\varepsilon\}.
\]
Then in each $\ID _\varepsilon$ we have
\[
\frac{\Xi_p}{|\Xi_p|}\to\frac{\psi}{|\psi|}
\]
locally uniformly. Since $\psi$ is holomorphic, we have $|\ID -\ID _\varepsilon|\to0$ as $\varepsilon\to0$. Then, as $\|Dh_p\|_{L^2(\ID )}$ are uniformly bounded, it follows that for any $\varepsilon>0$ and compact $A\subset\ID $,
\begin{align*}
&\lim_{p\to\infty}\int_{\ID _\varepsilon\cap A}\Big|(h_p)_\wbar \frac{\psi}{|\psi|}-k(h_p)_w\Big|\\
\leq&\lim_{p\to\infty}\int_{\ID _\varepsilon\cap A}\Big|(h_p)_\wbar \big(\frac{\psi}{|\psi|}-\frac{\Xi_p}{|\Xi_p|}\big)\Big|+\int_{\ID _\varepsilon\cap A}\Big|(h_p)_\wbar \frac{\Xi_p}{|\Xi_p|}-k(h_p)_w\Big|=0.
\end{align*}
Also,
\[
(h_p)_\wbar \frac{\psi}{|\psi|}-k(h_p)_w\rightharpoonup h_\wbar \frac{\psi}{|\psi|}-kh_w
\]
in $L^2(\ID )$. So we conclude
\[
h_\wbar \frac{\psi}{|\psi|}-kh_w=0.
\]
As $h$ is quasiconformal we have $|h_w|>0$ a.e. Thus
\[
\mu_h=k\frac{|\psi|}{\psi}=k\frac{\overline{\psi}}{|\psi|}
\]
almost everywhere in $\ID _\varepsilon\cap A$. By the arbitrariness of $\varepsilon$ and $A$ this holds almost everywhere in $\ID $.
\hfill $\Box$

\medskip

For nondegenerate Hamilton sequences the argument above,  which really follows Ahlfors, establishes a result usually attributed to Reich and Strebel stating that if $\mu_f$ has a Hamilton sequence that does not degenerate, then f is a Teichm\"uller mapping.  A proof of which can be found in \cite{EL}.
 
 \medskip
 
 We now have to establish Corollaries \ref{cor1} and \ref{cor2}.  First if $\phi_p\neq 0$,  then $ \phi_p/\|\phi_p\|_1\neq 0$ and so by Hurwitz's theorem the local uniform limit $\psi$ is either identically zero,  or nonzero.  The first case is ruled out by our assumption on nondegeneracy.  Thus $\mu_h=k\psi/|\psi|$ is $C^{\infty}(\ID)$ as $\psi\neq 0$.  This implies $h$ (and $f=h^{-1}$) are diffeomorphisms and thus  Corollary \ref{cor1} follows.  Next,  suppose $f_0$ are boundary values which do not admit a smooth extremal quasiconformal mapping (for instance if $f_0$ is not the boundary values of a M\"obius transformation,  but is   automorphic with respect to a co-compact Fuchsian group) but that the extremal quasiconformal mapping is a Teichm\"uller mapping with distortion $k\psi/|\psi|$.  Then unique extremality tells us that $h_p\to h$ and our argument above show $\phi_p\to\psi$. This sequence is non-degenerate with the choice of $f_0$ as described above since the space of such quadratic differentials is finite dimensional - see Ahlfors \cite[pp X]{Ahlfors}. Hurwitz's theorem tells us that $\phi_p$ must vanish for all sufficiently large $p$.   We observe from (\ref{7.2}) that $\phi_p(w_0)=0$ implies $h_\wbar=0$ as $\IK(w,h)\geq 1$ and $|h_w|\geq |h_\wbar|$. The local Lipschitz regularity of $h_p$,  established in \cite{IKO}, shows that $h_w\in L^{\infty}_{loc}(\ID)$ and hence $\mu_h=0$ implies $h_\wbar = 0$ and hence $\mu_h(w_0)=0$.  This establishes Corollary \ref{cor2}.

\subsection{$p\to 1$}

We presume that there are homeomorphisms $h_0:\IS\to\IS$ whose harmonic extension $P[h_0]$ to the disk has $\IK(z,P[h_0]^{-1})\in L^1(\ID)$ but not in any $L^p(\ID)$ for $p>1$.  Thus the exact assumptions needed on the boundary values for the result we seek to hold are unclear.  We thus suppose that $f_0:\IS\to\IS$ is quasisymmetric.  Then for each $p>1$, let $f_p:\ID\to\ID$ be the pseudo-inverse of the extremal map $h_p\in \mathsf{E}_p$ with boundary values $h_0=f_{0}^{-1}$. We show that $h_p:\ID\to\ID$ converges locally uniformly to the unique harmonic mapping with boundary values $h_0$ and also that $\mathsf{E}_p(h_p)\to \mathsf{E}_p(h_1)$.

\medskip

Let $H:\ID\to\ID$ be the Douady-Earle (or Ahlfors-Beurling) extension of $h_0$ to the disk, \cite{DE}.  Then $H:\ID\to\ID$ is quasiconformal, $H|\IS=h_0$ and $H$ is a real analytic diffeomorphism in $\ID$. Let $r<1$ and $\Omega_r=H^{-1}(\ID(r))$ and define  $H_r:\ID\to\ID$ as follows.  
\begin{equation}
H_r(z) = \left\{\begin{array}{ll} H(z), & z\in \ID\Omega_r; \\
P[{H|\partial\Omega_r}],  & z\in \Omega_r. \end{array}\right.
\end{equation}
Here $P[{H|\partial\Omega_r}]$ is the harmonic extension (or Poisson integral) of $H|\partial\Omega_r$.  Since the image $H(\Omega_r)$ is convex, $P[{H|\partial\Omega_r}]$ is a diffeomorphism \cite{Duren}.  Since $H|\partial\Omega_r$ is a real analytic diffeomorphism $P[{H|\partial\Omega_r}]$ is quasiconformal \cite{Martio}. We do not need any uniform bounds here,  simply that $H_r$ is now quasiconformal and so has an inverse with $p$-integrable distortion.  Now the sequence $\{h_p\}$ is uniformly bounded in $W^{1,2}(\ID)$, equicontinuous and therefore admits a subsequence converging locally uniformly to a mapping $h_1:\ID\to\ID$, $h_1\in W^{1,2}(\ID)$ and $h_1$ continuous as a map $\overline{\ID}\to\overline{\ID}$.  Strictly speaking we need a modulus of continuity estimate on a neighbourhood of $\overline{\ID}$.   This is achieved by reflection,  that is we consider 
\[ h_p^\ast = \left\{\begin{array}{ll} h_p(z), & z\in \overline{\ID}  \\
H(z/|z|^2)/|H(z/|z|^2)|^2  & z\in \IC \setminus \ID. \end{array}\right. \]
Then $h_p^\ast \in W^{1,2}(U)$ with uniform bounds for some open $U$ with  $\overline{\ID}\subset U$ as $H$ is quasiconformal. Polyconvexity of the integrand here gives
\begin{eqnarray*}
\int_\ID \|Dh_1\|^2  \;dw & \leq &\liminf_{p\to1}\int_\ID \| Dh_p \|^2 \;dw  = \liminf_{p\to1}\int_\ID \IK(w,h_p) J(w,h_p) \;dw   \\
& \leq & \liminf_{p\to1}\int_\ID \IK^p(w,h_p) J(w,h_p) \;dw  = \liminf_{p\to1} \mathsf{E}_p(h_p)
\end{eqnarray*}
where we have chosen to ignore the passage to a subsequence as it will not matter.  Notice also that by extremality  $\mathsf{E}_p(h_p)\leq \mathsf{E}_p(h_q) \leq \mathsf{E}_q(h_q) $ if $q\geq p$.
Next,
\begin{eqnarray*}
\lefteqn{\int_\ID \IK^p(w,h_p) J(w,h_p) \;dw}\\  & = & \mathsf{E}_p \leq   \int_\ID \IK^p(w,H_r) J(w,H_r) \;dw \\ 
& = & \int_{\ID\setminus\Omega_r} \IK^p(w,H) J(w,H)\; dw +\int_{\Omega_r} \IK^p(w,H_r) J(w,H_r)\; dw \\
& = & \int_{\ID\setminus\ID(r)} \IK^p(w,H^{-1})  \; dw + \int_{\Omega_r} \IK^p(w,H_r) J(w,H_r) \; dw \\
& = &\|\IK(z,H^{-1})^p\|_{\infty}(1-r^2) + \int_{\Omega_r} \IK^p(w,H_r) J(w,H_r) \; dw \\
\end{eqnarray*}
We now let $p\to 1$ in this inequality to obtain (with $M=\|\IK(z,H^{-1})\|_{\infty}$)  from the Monotone convergence theorem
\begin{equation} \liminf_{p\to1}\int_\ID \| Dh_p \|^2 \;dw \leq M(1-r^2) + \int_{\Omega_r} \IK(w,H_r) J(w,H_r) \; dw
\end{equation}
Now
\begin{eqnarray*}
\lefteqn{ \int_{\Omega_r} \IK(w,H_r) J(w,H_r) \; dw}\\
& = & \int_{\Omega_r}\|DH_r\|^2 \; dw = \int_{\Omega_r}\|D(P[H_r|\partial\Omega_r])\|^2 \; dw 
=  \int_{\Omega_r}\IK(z,(P[H_r|\partial\Omega_r])^{-1})\; dz \\ & \leq & \int_{\ID(r)} \IK(z,H^{-1})\; dz \leq  \int_{\ID} \IK(z,H^{-1})\; dz =M<\infty.
\end{eqnarray*}
Now, as $r\to 1$,  $\{H_r\}$ is a uniformly bounded sequence of harmonic mappings.  Their local uniform limit exists and is also harmonic and further, these energies converge as well as they are uniformly bounded. This limit is the Harmonic extension of $h_0$,  which we know to be diffeomorphic and a unique absolute minimiser \cite{AIMO}.  Putting together what we now have we see that as $r\to 1$ 
\begin{eqnarray*}
\mathsf{E}_1(h_1) \leq  \liminf_{p\to1} \mathsf{E}_p(h_p) \leq \mathsf{E}_1(P[h_0|\IS]) \leq \mathsf{E}_1(h_1)
\end{eqnarray*}
Thus equality holds throughout and $h_1=P[h_0|\IS] = \lim h_p$,  the convergence is uniform on $\ID$ and in $W^{1,2}(\ID)$ and 
\[ \mathsf{E}_p(h_p)\searrow \mathsf{E}_1(h_1)\]

\section{Examples.}

In this section we give a family of examples with non-constant Ahlfors-Hopf differential $\phi=\alpha \,z^{-2}$. When $m=1$ below,  we know these maps are uniquely extremal  in $F_p$ (for their boundary values) and are self-diffeomorphisms of $\ID$,  \cite{MJ}.

\begin{theorem} There is a smooth diffeomorphism $f_0:\IS\to\IS$ for which the minimiser $F_0\in F_p$ is a smooth quasiconformal diffeomorphism,  which is not conformal.
There is a homeomorphism $f_1:\IS\to\IS$ which is $ C^{\infty}$ smooth away from a single point for which  the minimiser $F_1 \in F_p$ is a diffeomorphism but it is not quasiconformal.  
\end{theorem}

Our examples will be found among the radial mappings defined on a sector.  Thus let
\begin{equation}
h(w)=H(r)e^{i m \theta}, \hskip10pt w=re^{i\theta},
\end{equation}
and $H:(r_0,R_0)\to(s_0,S_0)$ strictly increasing. We compute that
\begin{eqnarray*} 
h_w= \frac{1}{2} e^{i(m-1)\theta}\Big(\dot{H}+m \frac{H}{r}\Big), && 
h_\wbar= \frac{1}{2}  e^{i(m+1)\theta}\Big(\dot{H}-m \frac{H}{r} \Big)
\end{eqnarray*}
and so the
 Ahlfors-Hopf differential
\[ \IK^{p-1}(w,h)h_w \overline{h_\wbar} = \Big( \frac{r \dot{H}}{m H}+\frac{m H}{r \dot{H}}\Big)^{p-1} \Big(\dot{H}^2-\frac{m^2H^2}{r^2}\Big)e^{-2i\theta} \]
The only holomophic functions with argument $e^{-2i\theta}$ are of the form $\alpha z^{-2}$ and hence the Ahlfors-Hopf differential is holomorphic if and only if
\begin{equation}\label{8.2}
\Big( \frac{r \dot{H}}{mH}+\frac{mH}{r \dot{H}}\Big)^{p-1} \Big(r^2\dot{H}^2- m^2H^2 \Big)   = \alpha, 
\end{equation}
for some constant $\alpha$. If $\alpha=0$,  then $r\dot{H}=m H$ and $h(z)=cz^m$ is conformal.  

This equation takes a nicer form in terms of the inverse
\begin{equation}
f(z)=F(\rho) e^{i\theta}, \hskip10pt z=\rho e^{i\theta}
\end{equation}
and $F:(s_0,S_0)\to(r_0,R_0)$ strictly increasing.  Then (\ref{8.2}) becomes
\begin{equation}\label{8.3}
\rho^2\Big(1-\frac{mF^2 }{\rho^2\dot{F}^2 }\Big)\Big(\frac{\rho\dot{F}}{mF}+\frac{mF}{\rho\dot{F}}\Big)^{p-1}=\alpha,
\end{equation}
for $\alpha$ constant.  To study this equation we simplify with the assumption $m=1$ and define
\begin{equation}
a=a(\rho)=\frac{\rho \dot{F}}{F} > 0.
\end{equation}
Then (\ref{8.3}) reads as
\begin{equation}\label{8.4}
\rho^2(a^2-1 ) (a^2+1 )^{p-1}=\alpha a^{p+1},
\end{equation}
and this equation shows us that $0<a <1$ and $\alpha>0$,  or $1<a <\infty$ and $\alpha<0$.   It also defines $\rho$ as a smooth function of $a$, $a\neq1$,  increasing for $0<a<1$  and decreasing for $1<a<\infty$ and with range $(0,\infty)$ on both these intervals.
Since
$
(\log F)_\rho=\frac{\dot{F}}{F}=\frac{a}{\rho} 
$
we see that
\begin{eqnarray*}
\log F&=&\int\frac{a(\rho)}{\rho}d\rho=\int\frac{a}{\rho(a)}\rho_a(a)da =\int\frac{1+4a^2-a^4+(a^2-1)^2p}{2(1-a^4)}da\\
&=&-\frac{p-1}{2}a+(p-1)\arctan(a)+\frac{1}{2}\log\Big|\frac{a+1}{a-1}\Big| +C_1.
\end{eqnarray*}
for some constant $C_1$.  Let us choose the regime $\alpha<0$ and
\[ C_1= \frac{p-1}{2}a_1-(p-1)\arctan(a_1)-\frac{1}{2}\log\Big|\frac{a_1+1}{a_1-1}\Big|, \]
where $a_1<1$,  solves $(a_1^2-1 ) (a_1^2+1 )^{p-1}=\alpha a_1^{p+1}$ so that $F(1)=1$.  Then $f:\ID\setminus\{0\}\to A(s_0,1)$,  $s_0=e^{C_1}$.
Now $f$ is a quasiconformal diffeomorphism on any disk $D$ with $\overline{D}\subset \ID_\ast$. If $\phi:f(D)\to\ID$ is any Riemann map, and $\psi:\ID\to D$ is a similarity,   then $\phi\circ f\circ\psi :\ID\to\ID$ is uniquely extremal for its boundary values.  To see this if $g$ is another such map with the same or smaller energy,  then we can construct a new mapping $\tilde{f}:\ID\setminus\{0\}\to\IC$ by
\begin{equation}
\tilde{f} = \left\{ \begin{array}{ll} f(z); & z\not\in D, \\ \phi^{-1}\circ g\circ \psi^{-1}; &z\in D. \end{array}\right.
\end{equation}
It is straightforward to see that $\tilde{f}$ is a $W^{1,2}_{loc}(\ID)$ homeomorphism of $\ID\setminus\{0\}\to A(s_0,1)$.  Further,  if $D\subset A(r_0,1)$,  then $\tilde{f}:A(r_0,1)\to A(F(r_0),1)$ has energy no more than $f$.  However these radial maps are uniquely extremal \cite{MJ} and hence $\tilde{f}=f$ and $\phi\circ f\circ\psi =g$. \hfill $\Box$

\medskip

GJM Institute for Advanced Study, Massey University, Auckland, New Zealand \\
email: G.J.Martin@Massey.ac.nz

CY Institute for Advanced Study, Massey University, Auckland, New Zealand \\
email: C.Yao@massey.ac.nz 
\end{document}